                   \def\version{27 October, 2006}                            %

\documentclass[reqno,11pt]{amsart}
\usepackage{amsmath,epsfig}
\usepackage{amssymb}

\newtheoremstyle{rem}{1.3ex}{1.3ex}{\rmfamily}{}
{\itshape\rmfamily}{}{1.5ex}{}

\newenvironment{proofsect}[1]
{\vskip0.1cm\noindent{\bf #1.}}

\newtheorem{theorem}{Theorem}[section]
\newtheorem{lemma}[theorem]{Lemma}
\newtheorem{prop}[theorem] {Proposition}

\newtheorem{cor}[theorem] {Corollary}

\theoremstyle{definition}

\newtheorem{bem}[theorem] {Remark}

\renewcommand{\section}{\secdef\sct\sect}
\newcommand{\sct}[2][default]{\refstepcounter{section}
\setcounter{equation}{0}
\vspace{0.5cm}
\centerline{ \large
\scshape \arabic{section}.\ #1}
\vspace{0.3cm}}
\newcommand{\sect}[1]{
\vspace{0.5cm}
\centerline{\large\scshape #1}
\vspace{0.3cm}}

\renewcommand{\subsection}{\secdef \subsct\sbsect} 
\newcommand{\subsct}[2][default]{\refstepcounter{subsection} 
\nopagebreak 
\vspace{0.5\baselineskip} 
{\flushleft\bf \arabic{section}.\arabic{subsection}~\bf #1  } 
\nopagebreak} 
\newcommand{\sbsect}[1]{\vspace{0.1cm}\noindent 
{\bf #1}\vspace{0.1cm}} 

\newcommand{\Ocal}{{\mathcal {O}}}
\newcommand{\Dcal}{{\mathcal {D}}}
\newcommand{\Tcal}{{\mathcal {T}}}

\newcommand{\eps}{\varepsilon}
\newcommand{\id}{{\operatorname {id}}}

\newcommand{\R}     {\mathbb{R}}
\newcommand{\Z}     {\mathbb{Z}}
\newcommand{\N}     {\mathbb{N}}
\renewcommand{\P}   {\mathbb{P}}

\newcommand{\E}     {\mathbb{E}}

\renewcommand{\d}   {\operatorname{d}\!}

\newcommand{\sign}  {\operatorname{sign}}

\newcommand{\Sym}   {\mathfrak{S}}
\newcommand{\heap}[2]{\genfrac{}{}{0pt}{}{#1}{#2}}
\def\1{{\mathchoice {1\mskip-4mu\mathrm l}      
                    {1\mskip-4mu\mathrm l}
                    {1\mskip-4.5mu\mathrm l} {1\mskip-5mu\mathrm l}}}
\newcommand{\ssup}[1] {{\scriptscriptstyle{({#1}})}}

\begin{document}

\title[Ordered random walks]{\large Ordered random walks}

\author[Peter Eichelsbacher and Wolfgang K\"onig]{}
\maketitle
\thispagestyle{empty}
\vspace{0.2cm}

\centerline{\sc Peter Eichelsbacher\footnote{Ruhr-Universit\"at Bochum, Fakult\"at f\"ur Mathematik, NA3/68, D-44780 Bochum, Germany, {\tt Peich@math.ruhr-uni-bochum.de  }} and Wolfgang K\"onig\footnote{Universit\"at Leipzig, Mathematisches Institut, Augustusplatz 10/11, D-04109 Leipzig, Germany, {\tt koenig@math.uni-leipzig.de}\\ 
Both authors have been supported by Deutsche Forschungsgemeinschaft via SFB/TR 12.}}
\vspace{2 cm}

\centerline{\small{\version}}
\vspace{1cm}

\begin{quote}
{\small {\bf Abstract:} We construct the conditional version of $k$ independent and identically distributed random walks on $\R$ given that they stay in strict order at all times. This is a generalisation of so-called non-colliding or non-intersecting random walks, the discrete variant of Dyson's Brownian motions, which have been considered yet only for nearest-neighbor walks on the lattice. Our only assumptions are moment conditions on the steps and the validity of the local central limit theorem. The conditional process is constructed as a Doob $h$-transform with some positive regular function $V$ that is strongly related with the Vandermonde determinant and reduces to that function for simple random walk. Furthermore, we prove an invariance principle, i.e., a functional limit theorem towards Dyson's Brownian motions, the continuous analogue.}

\end{quote}

\vfill


\eject

\setcounter{section}{0}

\section{Introduction and main result}\label{intro}

\subsection{Dyson's Brownian motions and non-colliding processes.}\label{sec-DysonBM}

\noindent In 1962, F.~Dyson \cite{Dy62} made a beautiful observation. He looked at a process version of the famous Gaussian Unitary Ensemble (GUE), a matrix-valued diffusion known as Hermitian Brownian motion. He was interested in the process of the vectors of the eigenvalues of that matrix process. It turned out that this process admits a concise description: it is in distribution equal to a family of standard Brownian motions, conditional on having never any collision of the particles. More explicitly, it is the conditional distribution of $k$ independent standard Brownian motions $B_1,\dots,B_k$ on $\R$ given that the $k$-dimensional Brownian motion $B=(B_1,\dots,B_k)$ never leaves the {\it Weyl chamber},
\begin{equation}\label{Wdef}
W=\bigl\{x\in\R^k\colon x_1< x_2< x_3<\dots< x_k\bigr\}.
\end{equation}
In other words, $B$ is conditioned on the event $\{T=\infty\}$, where
\begin{equation}\label{Tdef}
T=\inf\{t\in[0,\infty)\colon B(t)\notin W\}
\end{equation}
is the first time of a collision of the particles. The definition of the conditional process needs some care, since the event $\{T=\infty\}$ has zero probability. As usual in such cases, it is defined via a Doob $h$-transform with some suitable harmonic function $h\colon W\to(0,\infty)$. It turned out that a suitable choice for $h$ (in fact, the only one, up to constant multiples)
is the {\it Vandermonde determinant\/} $\Delta\colon\R^k\to\R$ given by
\begin{equation}\label{Vdm}
\Delta(x)=\prod_{1\leq i<j\leq k}(x_j-x_i)=\det\bigl[(x_j^{i-1})_{i,j=1,\dots,k}\bigr],\qquad x=(x_1,\dots,x_k)\in\R^k.
\end{equation}
More precisely, $h=\Delta\colon W\to(0,\infty)$ is a positive harmonic function for the generator $\frac 12\sum_{i=1}^n\partial_i^2$ of $B=(B_1,\dots,B_k)$ on $W$, and $\Delta(B(t))$ is integrable for any $t>0$ under any starting measure of the motions. Hence we may consider the Doob $h$-transform of $B$ on $W$ with $h=\Delta$. The transformed process is called {\it Dyson's Brownian motions}. It is known since long that this transformed process is identical to the limiting conditional process given $\{T>t\}$ as $t\to\infty$. Therefore, the process is also called {\it non-colliding Brownian motions}.

For some decades after this discovery, it was quiet about non-colliding random processes, but the interest renewed in the 1990ies, and it has become an active research area and is being studied for a couple of reasons. M.-F.~Bru \cite{Br91} studied another interesting matrix-valued stochastic  process whose eigenvalue process admits a nice description in terms of non-colliding random processes, the {\it Wishart processes}, which are based on squared Bessel processes in place of Brownian motions. These processes and some few more were studied in \cite{KO01}. Non-colliding Brownian motions on the circle were investigated in \cite{HW96}, asymptotic questions about Brownian motions in a Weyl chamber in \cite{Gr00}, and a systematic study of a large class of physically relevant matrix-valued processes and their eigenvalue processes is carried out in \cite{KT04}. 

Certainly, also the time-discrete version has been studied, more precisely, families of $k$ i.i.d.~discrete random walks, conditional on never leaving $W$. It is important for the present paper to note that so far only random walks have been considered that have the following {\it continuity property}: at the first time of a violation of the strict ordering, there are two components of the walk that are at the same site (and produce therefore a collision). In other words, leaving $W$ is only possible via a step into the boundary $\partial W$ of $W$. This property is shared by nearest-neighbor random walks on the lattice $\Z^k$, started in $(2\Z)^k\cap W$ (in which case the walkers cannot jump over each other) and by walks that have only steps in $\{0,1\}^k$ or by imposing similar rules. Obviously, this continuity property makes the analysis much easier, but heavily restricts the choice of the step distribution. For walks having this property, the event of never leaving $W$ (i.e., of being strictly ordered at any time) is identical to being non-colliding, hence the term {\it non-colliding random walks\/} became popular, but also {\it vicious walkers}, {\it non-intersecting paths} and {\it non-intersecting walks}. We consider the latter two terms misleading since it is the {\it graphs\/} that are non-intersecting, more precisely the graph of the polygon line that interpolates between discrete time units. Non-intersecting paths played an important role in the proof of Johanssons's beautiful analysis \cite{Jo00}, \cite{Jo02} of the corner-growth model (which is equivalent to directed first-passage percolation). These works naturally raise the interesting question how far the connections between the corner-growth model and non-intersecting paths reach; they are yet known only for rather restricted waiting-time distributions respectively step distributions. Further relationships to other models, like the Arctic circle, are investigated in \cite{Jo02}. Recently \cite{BS06}, a random matrix central limit behavior was obtained for the rescaled versions of many non-intersecting random walks with an essentially general step distribution. The non-intersecting property was required only up to a fixed time. Furthermore, also applications in the study of series of queues in tandem were found and analysed; see the survey article \cite{OC03}. 

Especially in recent years, more and more connections have been found between non-colliding random processes and various models, some of which have not yet been fully understood. A number of explicit examples have been worked out, and the class of random processes whose non-colliding version could be rigorously established and characterized, is growing. It is now known how to construct and describe these conditional versions for a couple of examples of random walks, among which the binomial random walk, the multinomial walk, and the (continuous-time) Poisson random walk \cite{KOR02}, and birth and death processes and the Yule process \cite[Ch.~6]{Do05}. In all these explicit examples, it fortunately turned out that the Vandermonde determinant, $\Delta$, is a positive regular function for the generator of the family of the random walks, and the Doob $h$-transform with $h=\Delta$ could explicitly be calculated. A survey on non-colliding random walks appears in \cite[Ch.~4]{K05}.

However, to the best of our knowledge, the theory of non-colliding random processes still consists of a list of explicit, instructive and important examples, but the general picture is still lacking. In particular, the precise class of random walks for whose generator the Vandermonde determinant is  a positive regular function, is widely unknown yet, and it is also yet unknown what function in general replaces $\Delta$ in the construction, if it can be carried out. 

The present paper reveals the general mechanism of constructing from a tuple of $k$ i.i.d.~random walks on $\R$ the conditional version that never leaves $W$, i.e., whose components stay in strict order at any time. Only the finiteness of some sufficiently high  moments of the walker's steps and the validity of the local central limit theorem will be assumed. We will identify a positive harmonic function in terms of which we will construct the version that never leaves $W$. Furthermore, we will also consider the asymptotic behavior of the conditional walk and prove an invariance principle, i.e., the convergence of the properly rescaled process towards the continuous version, Dyson's Brownian motions. We consider the results of this paper as a universal approach to non-intersecting paths, which opens up a possibility to attack in future also related models like the corner-growth model in a universal manner.

Since a general random walk makes jumps of various sizes, the term {\it non-colliding} is misleading, and the term {\it non-intersecting} refers to the graphs instead of the walks. We prefer to replace these terms by {\it ordered random walks}, for obvious reasons. General non-colliding random walks in the strict sense seem to represent an exciting and open research topic that may be inspired from other topics than processes of random matrices and will presumably not have much to do with Dyson's Brownian motions.

\subsection{Ordered random walks}\label{sec-ordered}

\noindent For $k\in\N$, let $X_1,\dots,X_k$ be $k$ independent copies of a random walk, $X_i=(X_i(n))_{n\in\N_0}$, on $\R$. Then $X=(X_1,\dots,X_k)$ is a random walk on $\R^k$ with i.i.d.~components. Our goal is to construct a conditional version of $X$, given that the $k$ components stay in a fixed order for all times. That is, we want to condition $X$ on never leaving the  Weyl chamber $W$ in \eqref{Wdef}. Another way to formulate this is to condition on the event $\{\tau=\infty\}$, where
\begin{equation}
\tau=\inf\{n\in\N_0\colon X(n)\notin W\}
\end{equation}
is the first time that some component reaches or overtakes another one.  Some care is needed in defining the conditional process, since the event $\{\tau=\infty\}$ has zero probability. We shall construct this process as a Doob $h$-transform and show that it coincides with the limiting conditional process given $\{\tau>n\}$ as $n\to\infty$.

Let $S\subset\R$ denote the state space of the random walk $X_1$ when started at 0. Let $\P$ denote the underlying probability measure. For $x\in\R^k$, we write $\P_x$ when the process $X=(X(n))_{n\in\N_0}$ starts at $X(0)=x$, and we denote by $\E_x$ the corresponding expectation. A function $h\colon W\cap S^k\to(0,\infty)$ is called a {\it positive regular function\/} with respect to the restriction of the transition kernel of $X$ to $W\cap S^k$ if
\begin{equation}
\E_x[h(X(1))\1_{\{\tau>1\}}]=h(x),\qquad x\in W\cap S^k.
\end{equation}
In this case, we may define the Doob $h$-transform of the process $X$ via the $n$-step transition probabilities
\begin{equation}\label{trans}
\widehat \P^{\ssup h}_x(X(n)\in\d y)=\P_x(\tau>n;X(n)\in\d y)\frac {h(y)}{h(x)},\qquad x,y\in W\cap S^k, n\in\N.
\end{equation}
The regularity and positivity of $h$ guarantee that the right hand side of \eqref{trans} is a probability measure on $W$ in $\d y$. The state space of the Doob 
$h$-transform is equal to $W\cap (S^k-x)$ when started at $x$. {\it A priori}, the existence and uniqueness of such positive regular function is far from clear,
and also the question if the corresponding Doob transform has anything to do with the conditional version given $\{\tau >n\}$ in the limit as $n \to \infty$.

In the present paper, we present a positive regular function $V$ such that the Doob $h$-transform with $h=V$ turns out to be the conditional version of $X$ given never exiting $W$. Under the latter process, we understand (in the case of its existence) the limiting process $X$ given $\{\tau>n\}$ as $n\to\infty$. Furthermore, we analyse the decay of the probability of the event $\{\tau>n\}$ and give a limit theorem for the rescaled path's endpoint, $n^{-1/2}X(n)$, conditioned on this event. Another main goal is the analysis of the conditional process at large times. We show that the rescaled conditional process $(n^{-1/2}X(\lfloor t n\rfloor))_{t\geq 0}$ converges towards Dyson's Brownian motions.

Now we state the precise assumptions on the walk. We want to work with a random walk that lies in the domain of attraction of Brownian motion. Without loss of generality we therefore put the 

\medskip

\noindent {\bf Centering Assumption. } {\it The walk's steps have mean zero and variance one.}
 
\medskip
 
We distinguish the two cases of a lattice walk and a non-lattice walk. The following assumption will enable us to apply a local central limit theorem, which will be an important tool in our proofs.

\medskip

\noindent{\bf Regularity Assumption. } {\it Either the support of the walk, $\{X_1(n)\colon n\in\N_0\}$, under $\P_0$, is contained in the lattice $\alpha \Z$ for some minimal $\alpha\in(0,\infty)$, or the distribution of $X_1(N)$ possesses a bounded density for some $N\in\N$.}

\medskip

The walk's state space, $S$, is equal to $\alpha \Z$ in the first case, the lattice case, and it is equal to $\R$ in the second case, the non-lattice case.

Now we introduce the main object of the paper, the positive regular function $h=V$ we will be working with. Define $V\colon W\cap S^k\to\R$ by
\begin{equation}\label{Vdef}
V(x)=\Delta(x)-\E_x\bigl[\Delta(X(\tau))\bigr],\qquad x\in W\cap S^k.
\end{equation}
Actually, it is {\it a priori\/} not clear at all under what assumptions $V$ is well-defined, i.e., under what assumptions 
$\Delta(X(\tau))$ is integrable under $\P_x$ for any $x\in W$. This question is trivially answered in the affirmative for walks 
that have the above mentioned continuity property, which may be also formulated by saying that $\P_x(X(\tau)\in \partial W)=1$. 
This property depends on the initial site $x\in W$ (e.g.~simple random walk in $\Z^k$ starting in $(2\Z)^k$ has this property, but not 
when it starts in the site $(1,2,\dots,k)$, say). 
All examples of walks considered in the literature so far (see Section~\ref{sec-DysonBM}) have this property. If the walk has this property, then $X(\tau)$ has some equal components, $\P_x$-a.s., and therefore the Vandermonde determinant $\Delta(X(\tau))$ equals zero, $\P_x$-a.s., which shows that $V(x)=\Delta(x)$. 

However, in the general case considered in the present paper, the integrability of $\Delta(X(\tau))$ seems subtle, and we succeeded in proving the integrability only under some moment condition on the steps and the local central limit theorem.

\begin{theorem}\label{main} Assume that the random walk $X$ satisfies the Centering Assumption and the Regularity Assumption. Then, there is a $\mu=\mu_k>0$, depending only on $k$, such that, if the $\mu$-th moment of the walk's steps is finite, the following hold.
\begin{enumerate}
\item[(i)] For any $x\in W$, the random variable $\Delta(X(\tau))$ is integrable under $\P_x$.

\item[(ii)] The function $V$ defined in \eqref{Vdef} is a positive regular function with respect to the restriction of the transition kernel to $W\cap S^k$, and $V(X(n))$ is integrable with respect to $\P_x$ for any $x\in W$ and any $n\in\N$. 

\item[(iii)] The Doob $h$-transform of $X$ on $W\cap S^k$ with $h=V$ is equal to the distributional limit of the conditional process given $\{\tau>n\}$ as $n\to\infty$.

\item[(iv)] For any $x\in W$, the distribution of $n^{-\frac 12} X(n)$ under $\P_x(\,\cdot\mid \tau>n)$ converges towards the distribution on $W$ with density $y\mapsto\frac 1{Z_1}{\rm e}^{-\frac 12|y|^2}\Delta(y)$ (with $Z_1$ the norming constant). Moreover,
\begin{equation}
\lim_{n\to\infty}n^{\frac k 4(k-1)}\P_x(\tau>n)= K V(x),\qquad\mbox{where }K=\int_W\frac{{\rm e}^{-\frac 12|y|^2}}{(2\pi)^{k/2}}\Delta(y)\,\d y\,\prod_{l=0}^{k-1}\frac 1{l!}.
\end{equation}

\item[(v)] For any  $M>0$, uniformly in  $x\in W$ satisfying $|x|\leq M$, $\lim_{n\to\infty}n^{-\frac k4(k-1)}V(\sqrt n \,x)=\Delta(x)$.

\item[(vi)] For any $x\in W$, the distribution of $n^{-\frac 12} X(n)$ under $\widehat \P^{\ssup {V}}_{x}$ converges towards the distribution on $W$ with density $y\mapsto\frac 1{Z_2}{\rm e}^{-\frac 12|y|^2}\Delta(y)^2$, the {\em Hermite ensemble} (with $Z_2$ the norming constant). More generally, the distribution of the process $(n^{-\frac 12} X(\lfloor nt\rfloor)_{t\in[0,\infty)}$ under the transformed probabilities, i.e., under $\widehat \P^{\ssup {V}}_{\sqrt n x}$ for any $x\in W$, converges towards Dyson's Brownian motions started at $x$.

\end{enumerate}
\end{theorem}

The proof of Theorem~\ref{main} is distributed over a couple of propositions and lemmas. More precisely, (i) is contained in Proposition~\ref{Xtaukint}, (ii) in Lemma~\ref{Vbound}, (iii) in Lemma~\ref{lem-CondInt}, (iv) in Corollary~\ref{cor-tautail}, (v) in Lemma~\ref{lem-VasyDelta} and (vi) in Lemma~\ref{lem-Dysonconv}. An explicit formula for the transition probabilities of the transformed process appears in \eqref{transfProc} below.

The only role of the Regularity Assumption is to establish the expansion in the local central limit theorem in \eqref{locexp} below, under a sufficient moment condition. Hence, all conclusions of Theorem~\ref{main} hold under \eqref{locexp} instead of the Regularity Assumption.

Our main tool in the proof of Theorem~\ref{main} is an extension of the well-known Karlin-McGregor formula to arbitrary random walks on $\R^k$ with i.i.d.~step distributions. Furthermore, we use H\"older's inequality (this is why we lose control on the minimal integrability assumption), the local central limit theorem and Donsker's invariance principle. Another helpful fact, proved in \cite{KOR02}, is that the process $(\Delta(X(n)))_{n\in\N}$ is a martingale, provided that $\Delta(X(n))$ is integrable for any $n\in\N$. 

Future work will be devoted to the study of the system of $k$ ordered walks in the limit $k\to\infty$, from which we hope to deduce interesting and universal variants of Wigner's semicircle law. The problem remains open under what minimal assumptions the assertions of Theorem~\ref{main} remain true. It seems that the integrability of $\Delta(X(\tau))$ requires at least the finiteness of the $k$-th moments of the steps. The special case $k=2$ includes the well-known and much-studied question of conditioning a single path to be positive at all times. Here it is known that the finiteness of the first moment of the step distribution is sufficient for the existence of a positive regular function. Furthermore, the finiteness of the $(l+1)$-st moment of the step distribution implies the finiteness of the $l$-th moment of $\Delta(X(\tau))$, for any $l\in\N$. The standard proof of these facts uses fluctuation theory and the Sparre-Andersen identity, see \cite[Chapter XII and XVIII]{F71}, e.g. Instead, we use rough moment estimates and H\"older's inequality in the present paper and therefore lose control on minimal integrability assumptions. To the best of our knowledge, it is open how to construct the ordered version in the case of even less integrability.

The remainder of this paper is devoted to the proof of Theorem~\ref{main}. In Section~\ref{sec-KMG}, we present our main tool, a generalisation of the well-known Karlin-McGregor formula, a determinantal formula for the marginal distribution before the first time of a violation of the strict ordering. In Section~\ref{sec-Vexist} we prove that $\Delta(X(\tau))$ is integrable under $\P_x$ for any $x\in W$, a fact which establishes that $V(x)$ is well-defined. Finally, in Section~\ref{sec-Vpos} we prove a couple of properties of $V$, in particular its positivity (a fact which is crucial to define the transformed process) and the functional limit theorem towards Dyson's Brownian motions.

\section{A generalized Karlin-McGregor formula}\label{sec-KMG}

An important tool for handling the distribution of the process $X$ on the event $\{\tau>n\}$ is the well-known {\it Karlin-McGregor formula\/} \cite{KM59} for the transition probabilities before a violation of the strict ordering of the components. This is an explicit formula for the distribution of $X(n)$ on $\{\tau>n\}$ for a variety of stochastic processes including nearest-neighbor random walks on $\Z^k$ and Brownian motion. In the standard Brownian motion case, this formula reads
\begin{equation}\label{KMold}
{\tt P}_x(T>t;B(t)\in\d y)=\det\Big[\Big({\tt P}_{x_i}(B_1(t)\in\d y_j)\Big)_{i,j=1,\dots,k}\Big],\qquad t>0,x,y\in W,
\end{equation}
where the $k$ motions start from $x$ under ${\tt P}_x$ (recall \eqref{Wdef} and \eqref{Tdef}).
The proof is based on the continuity of the paths and on the reflection principle: if the two motions $B_i$ and $B_j$ meet each other at time $T=s\in(0,t)$, then the paths $ (B_j(r))_{r\in[s,t]}$ and $ (B_i(r))_{r\in[s,t]}$ are interchanged, and we obtain motions that arrive at $y_j$ and $y_i$ rather than at $y_i$ and $y_j$. A clever enumeration shows that \eqref{KMold} holds. For this method to work it is crucial that the two motions $B_i$ and $B_j$ are located at the same site at time $T$. The same argument applies to many other processes including discrete-time walks on the lattice $\Z^k$ that have the continuity property discussed prior to Theorem~\ref{main}, i.e., $\P_x(X(\tau)\in\partial W)=1$. Since the proof of the Karlin-McGregor formula also involves a reflection argument, it is valid only for walks on $\Z^k$ whose step distribution is i.i.d. 

In the present paper, we overcome the continuity restriction, and we work with ($k$ copies of) an arbitrary walk on the real line. 
We succeed in finding an analogue of the formula in \eqref{KMold}, which we present now. Introduce the signed measure
\begin{equation}\label{Dndef}
\begin{aligned}
{\Dcal}_n(x,\d y)&=\det\Big[\Big(\P_{x_i}(X_1(n)\in\d y_j)\Big)_{i,j=1,\dots,k}\Big]\\
&=\sum_{\sigma\in\Sym_k}\sign(\sigma) \prod_{i=1}^k\P_{x_{\sigma(i)}}(X_1(n)\in\d y_i),\qquad x=(x_1,\dots,x_k), \, y=(y_1,\dots,y_k),
\end{aligned}
\end{equation}
where $\Sym_k$ denotes the set of permutations of $1,2,\dots,k$. 

The following is a generalization of \eqref{KMold} to general random walks on the real line. We use the notation of Section~\ref{sec-ordered}, but no assumptions on drift or moments are made, nor on existence of densities. 

\begin{prop}[Generalized Karlin-McGregor formula]\label{KMG} Let $(X(n))_{n\in\N_0}$ be an arbitrary random walk on $\R^k$ with i.i.d.~components. Then the following hold.
\begin{enumerate}
\item[(i)] For any $n\in\N$ and any $x,y\in W$,
\begin{equation}\label{KMGneu}
\P_x(\tau>n;X(n)\in\d y)={\Dcal}_n(x,\d y)-\E_x\bigl[\1_{\{\tau\leq n\}}{\Dcal}_{n-\tau}(X(\tau),\d y)\bigr].
\end{equation}

\item[(ii)] Define $\psi\colon \R^k\setminus W\to\R^k$ by 
\begin{equation}\label{psidef}
\psi(y)=(y_j-y_i)({\rm e}_j-{\rm e}_i),\qquad \mbox{where }(i,j)\in \{1,\dots,k\}^2 \mbox{ minimal  satisfying }i<j\mbox{ and } y_i>y_j.
\end{equation}
Here ${\rm e}_i$ is the $i$-th canonical unit vector in $\R^k$, and \lq minimal\rq\ refers to alphabetical ordering. Then, for any $l,n\in\N$ satisfying $l\leq n$ and any $x,y\in W$,
\begin{equation}\label{KMGneuneu}
-\E_x\bigl[\1_{\{\tau=l\}}{\Dcal}_{n-\tau}(X(\tau),\d y)\bigr]=
\E_x\Bigl[\1_{\{\tau=l\}}{\Dcal}_{n-\tau}\big(X(\tau),\d\; (y+\psi(X(\tau)))\big)
\Bigr].
\end{equation}
\end{enumerate}
\end{prop}

It is the assertion in (i) which we will be using in the present paper; no reflection argument is involved. The assertion in (ii) uses the reflection argument and is stated for completeness only.

In the special case of walks on $\Z$ that enjoy the above mentioned continuity property, $\P_x(X(\tau)\in\partial W)=1$, the second term on the right of \eqref{KMGneu} vanishes identically since the vector $X(\tau)$ has two identical components, and therefore the determinant vanishes.

An extension of Proposition~\ref{KMG} may be formulated for arbitrary Markov chains on $\R^k$ that satisfy the strong Markov property; the assertion in (ii) additionally needs exchangeability of the step distribution. However, the analogue of ${\Dcal}_n$ does not in general admit a {\it determinantal\/} representation.

\begin{proofsect}{Proof of Proposition~\ref{KMG}} We write $y_\sigma=(y_{\sigma(1)},\dots,y_{\sigma(k)})$. Using \eqref{Dndef}, we have
\begin{equation}\label{KMG1}
\begin{aligned}
\P_x(\tau>n;& \, X(n)\in \d y)-\Dcal_n(x,\d y)\\
&=\sum_{\sigma\in\Sym_k}\sign(\sigma)\bigl[\P_x(\tau>n;X(n)\in \d y_\sigma)-\P_x(X(n)\in \d y_\sigma)\bigr]\\
&=-\sum_{\sigma\in\Sym_k}\sign(\sigma)\P_x(\tau\leq n;X(n)\in\d y_\sigma),
\end{aligned}
\end{equation}
since all the summands $\P_x(\tau>n;X(n)\in\d y_\sigma)$ are equal to zero with the exception of the one for $\sigma=\id$. Apply the strong Markov property to the summands on the right hand side of \eqref{KMG1} at time $\tau$, to obtain
\begin{equation}
\P_x(\tau\leq n;X(n)\in\d y_\sigma)= \sum_{m=1}^n\int_{\R^k\setminus W}\P_x(\tau=m; X(m)\in\d z)\P_z(X(n-m)\in\d y_\sigma).
\end{equation}
Substitute this in \eqref{KMG1}, we obtain
\begin{equation}\label{KMG2}
\begin{aligned}
{ }\mbox{right side of \eqref{KMG1}}& =-\sum_{m=1}^n\int_{\R^k\setminus W} \P_x(\tau=m; X(m)\in \d z)\sum_{\sigma\in\Sym_k}\sign(\sigma)\P_z(X(n-m)\in \d y_{\sigma})\\
&=-\sum_{m=1}^n\int_{\R^k\setminus W} \P_x(\tau=m; X(m)\in\d z)\det\Bigl[\Big(\P_{z_i}(X_1(n-m)\in \d y_j)\Big)_{i,j=1,\dots,k}\Bigr]\\
&=-\sum_{m=1}^n \E_x\big[\1_{\{\tau=m\}} \Dcal_{n-m}(X(m),\d y)\big].
\end{aligned}
\end{equation}
This shows that (i) holds. In order to show (ii), we use the reflection argument of \cite{KM59}. Fix $l\in\{0,1,\dots,n\}$ and a continuous bounded function $f\colon\R^k\to\R$, then it is sufficient to show
\begin{equation}\label{reflection}
-\sum_{\sigma\in\Sym_k}\sign(\sigma)\E_x\Big[\1_{\{\tau=l\}}\E_{X(l)}\big[f(X_\sigma(n-l))\big]\Big]
=\sum_{\sigma\in\Sym_k}\sign(\sigma)\E_x\Big[\1_{\{\tau=l\}}\E_{y}\big[f(X_\sigma(n-l)+\psi(y))\big]\big|_{y=X(l)}\Big].
\end{equation}
This is done as follows. Given a transposition $\lambda=(i,j)$ satisfying $i<j$, let $\tau_\lambda=\inf\{n\in\N\colon X_i(n)\geq X_j(n)\}$ be the first time at which the $i$-th and the $j$-th component of the walk are not in strict order anymore. On the event $\{\tau=l\}$, there is a minimal transposition $\lambda$ such that $\tau=\tau_\lambda$. On the event $\{\tau=\tau_\lambda=l\}$, abbreviating $y=X(l)$, in the inner expectation we reflect the path $(y=X_\sigma(0), X_\sigma(1),\dots,X_\sigma(n-l))$ in the $i$-$j$-plane around the main diagonal (i.e., we interchange all the steps of the $i$-th component with the ones of the $j$-th component), and we obtain a path that terminates after $n-l$ steps at $X_{\sigma\circ\lambda}(n-l)+\psi(y)$. Obviously, the reflection is measure-preserving, and therefore we have
$$
\E_{y}\big[f(X_\sigma(n-l))\big]=\E_{y}\big[f(X_{\sigma\circ\lambda}(n-l)+\psi(y))\big],\qquad \mbox{a.s.\ on }\{\tau=\tau_\lambda=l\}\mbox{, where }y=X(l).
$$
Hence, denoting the set of transpositions by $\Tcal_k$, we have
$$
\begin{aligned}
-\sum_{\sigma\in\Sym_k}&\sign(\sigma)\E_x\Big[\1_{\{\tau=l\}}\E_{X(l)} \big[f(X_\sigma(n-l))\big]\Big]\\
&=\sum_{\lambda\in\Tcal_k}\sum_{\sigma\in\Sym_k}\sign(\sigma\circ\lambda)\E_x\Big[\1_{\{\tau=\tau_\lambda=l\}}\E_{y}\big[f(X_{\sigma\circ\lambda}(n-l)+\psi(y))\big]\big|_{y=X(l)}\Big].
\end{aligned}
$$
Now substitute $\sigma\circ\lambda$, interchange the two sums and carry out the sum on $\lambda$, to see that the right hand side is equal to the right hand side of \eqref{reflection}.
\end{proofsect}
\qed

\section{Existence of ${{V(x)}}$}\label{sec-Vexist}

In this section, we assume that $(X(n))_{n\in\N_0}$ is a random walk satisfying the assumptions of Theorem~\ref{main}. Furthermore, we fix $x\in W$. We prove that $V(x)$ in \eqref{Vdef} is well-defined. This is equivalent to showing the integrability of $\Delta(X(\tau))$ under $\P_x$. This turns out to be technically nasty and to require a couple of careful estimates. 
The proof of the integrability of $\Delta(X(\tau))$ will be split into a number of lemmas. 
In Section~\ref{sec-tails} we explain the subtlety of the problem and reduce the proof of the integrability of $\Delta(X(\tau))$ to the control of the tails of $\tau$. In Section~\ref{sec-expansions} we provide a version of a higher-order local central limit theorem for later use. Our main strategy is explained in Section~\ref{sec-strategy}, where we also formulate and prove the main steps of the proof. Finally, in Section~\ref{sec-integrability} we finish the proof of the integrability of $\Delta(X(\tau))$.

\subsection{Integrability of $\boldsymbol{\Delta(X(\tau))}$ and the tails of $\boldsymbol\tau$.}\label{sec-tails}
 
\noindent The reason that the proof of integrability of $\Delta(X(\tau))$ is subtle comes from the following heuristic observation. We want to show that the series
\begin{equation}\label{heuristik}
\sum_{n\in\N}\E_x[|\Delta(X(\tau))|\1_{\{\tau=n\}}]
\end{equation}
converges. Since we shall prove that $\P_x(\tau>n)\asymp n^{-\frac k4(k-1)}$, one can conjecture (but we actually do not prove that) that the event $\{\tau=n\}$ should have probability of order $n^{-\frac k4(k-1)-1}$. On this event, $|X(\tau)|$ is of order $|X(n)|\approx \sqrt n$, according to the central limit theorem. Therefore also all the differences $|X_i(n)-X_j(n)|$ with $1\leq i<j\leq k$ should be of that order, with one exception: at time $\tau=n$, there is a random pair of indices $i^*,j^*$ such that $X_{i^*}(n)$ and $X_{j^*}(n)$ are close together, since the $i^*$-th and $j^*$-th walk just crossed each other. Hence, $|\Delta(X(\tau))|$ should be of order $n^{-\frac 12}\prod_{1\leq i<j\leq k}\sqrt n=n^{\frac k4(k-1)-\frac 12}$, where the first term accounts for that random pair $i^*,j^*$. Hence, in the expectation in \eqref{heuristik}, there is a subtle extinction between two terms. However, that term should be of order $n^{-\frac 32}$ and therefore summable.

The next lemma shows that, for proving the finiteness of the series in \eqref{heuristik}, it will be crucial to control the tails of $\tau$.

\begin{lemma}\label{integr} Assume that the $\mu$-th moment of the steps is finite, for some $\mu>(k-1)(\frac k2(k-1)+2)$. Then there are $r\in(0,\frac k4(k-1)-1)$ and $\lambda\in(0,1)$ such that, for any set $M\subset W$ that is bounded away from zero, there is $C>0$ such that
\begin{equation}\label{summableneu}
\E_x\big[|\Delta(X(\tau))|\1_{\{\tau\leq n\}}\big]\leq C |x|^{(1+a)\frac k2(k-1)}+
C \Big(\sum_{l=\lceil |x|^{2(1+a)}\rceil}^n l^r \P_x(\tau >l)\Big)^\lambda ,\qquad n\in\N, x\in M, a\geq 0.
\end{equation}
\end{lemma}

Using \eqref{summableneu} for $a=0$ and using some obvious estimates, we also have that for any compact set $M\subset W$, there is $C>0$ such that
\begin{equation}\label{summable}
\E_x\big[|\Delta(X(\tau))|\1_{\{\tau\leq n\}}\big]\leq C+
C \Big(\sum_{l=1}^n l^r \P_x(\tau >l)\Big)^\lambda ,\qquad n\in\N, x\in M.
\end{equation}
We will use \eqref{summable} later in the present section and \eqref{summableneu} turns out to be crucial in the proof of Lemma~\ref{lem-VasyDelta} below.
\begin{proofsect}{Proof} For $l\in\N$, let $Y_i(l)$ be the $l$-th step of the $i$-th walk, hence, $X_i(n)=x_i+\sum_{l=1}^n Y_i(l)$, $\P_x$-almost surely, and all the random variables $Y_i(l)$ with $i\in\{1,\dots,k\}$ and $l\in\N$ are i.i.d.~with $\E[|Y_i(l)|^\mu]<\infty$.

We split the expectation on the left hand side of \eqref{summableneu} into the events $\{\tau\leq |x|^{2(1+a)}\}$ and $\{|x|^{2(1+a)}<\tau\leq n\}$. On the first event, we basically use that $|X(\tau)|\leq \Ocal(|x|^{1+a})$. Indeed, abbreviating $S_i(l)=Y_i(1)+\dots+Y_i(l)$, 
we obtain
$$
\begin{aligned}
\E_x\big[|\Delta(X(\tau))|\1_{\{\tau\leq |x|^{2(1+a)}\}}\big]
&\leq \sum_{l\leq |x|^{2(1+a)}}\E_x\Big[\prod_{i<j}\big|x_i-x_j+S_i(l)-S_j(l)\big|\1_{\{\tau=l\}}\Big]\\
&\leq \E_x\Big[\prod_{i<j}\Big(|x_i-x_j|+\max_{l\leq |x|^{2(1+a)}}|S_i(l)-S_j(l)|\Big)\Big].
\end{aligned}
$$
Under our moment assumption, the expectation of $\prod_{(i,j)\in A}\max_{l\leq |x|^{2(1+a)}}|S_i(l)-S_j(l)|$ is at most of order $|x|^{|A|(1+a)}$, for any $A\subset\{(i,j)\in \{1,\dots,k\}^2\colon i<j\}$. This explains the first term on the right hand side of \eqref{summableneu}.

Now we turn to the expectation on the event $\{|x|^{2(1+a)}<\tau\leq n\}$.
Define, for $m=1,\dots,k-1$,
\begin{equation}\label{tauijdef}
\tau_{m}=\inf\{n\in\N_0\colon X_m(n)\geq X_{m+1}(n)\},
\end{equation}
the first time at which the $m$-th and the $(m+1)$-st component violate the strict ordering.
Hence, $\tau=\inf_{m=1,\dots, k-1}\tau_{m}$. On $\{\tau =\tau_{m}=l\}$, we have
\begin{equation}\label{overshoot}
0\leq X_m(\tau_{m})-X_{m+1}(\tau_{m})\leq Y_{m}(l)-Y_{m+1}(l).
\end{equation}
We use H\"older's inequality twice as follows. Fix $p,q>1$ satisfying $\frac 1p+\frac 1q=1$. For any $\xi>0$, we have (abbreviating $\overline Y_{m,l}=Y_{m}(l)-Y_{m+1}(l)$),
\begin{equation}\label{Hoeldertrick}
\begin{aligned}
\E_x&[|\Delta(X(\tau))|\1_{\{|x|^{2(1+a)}<\tau=\tau_m\leq n\}}]\\
&\leq \E_x\Big[\sum_{l=\lceil|x|^{2(1+a)}\rceil}^n\Big(l^{(\frac k2(k-1)-1)(\frac 12 +\xi)} \1_{\{\tau=l\}}\Big)\Big(|Y_{m+1}(l)-Y_m(l)|\prod_{(i,j)\not=(m,m+1)}\frac{|X_i(l)-X_j(l)|}{l^{\frac 12+\xi}}\Big)\Big]\\
&\leq \E_x\Big[\Big(\sum_{l=\lceil|x|^{2(1+a)}\rceil}^n \big(l^{(\frac k2(k-1)-1)(\frac 12 +\xi)} \1_{\{\tau=l\}}\big)^p\Big)^{1/p} \\
&\qquad\qquad\times\Big(\sum_{l=\lceil|x|^{2(1+a)}\rceil}^n |\overline Y_{m,l}|^q \prod_{(i,j)\not=(m,m+1)}\Big|\frac{X_i(l)-X_j(l)}{l^{\frac 12+\xi}}\Big|^q\Big)^{1/q}\Big]\\
&\leq \E_x\Big[\tau^{p(\frac k2(k-1)-1)(\frac 12 +\xi)}\1_{\{|x|^{2(1+a)}<\tau\leq n\}}\Big]^{1/p} \\
&\qquad\qquad\times
\E_x\Big[\sum_{l=\lceil|x|^{2(1+a)}\rceil}^n l^{-q\xi[\frac k2(k-1)-1]}|\overline Y_{m,l}|^q\prod_{(i,j)\not=(m,m+1)}\Big|\frac{X_i(l)-X_j(l)}{\sqrt l}\Big|^q\Big]^{1/q}.
\end{aligned}
\end{equation}
We put $r=p(\frac k2(k-1)-1)(\frac 12 +\xi)-1$ and $\lambda=1/p$ and choose $\xi$ so small and $p$ so close to one that $r<\frac k4(k-1)-1$ and $q\xi[\frac k2(k-1)-1]>1$. When we pick $\xi=(\frac k2(k-1)-1)^{-1}(\frac k2(k-1)+2)^{-1}$, then this is achieved by any choice of $q>(\xi[\frac k2(k-1)-1])^{-1}=\frac k2(k-1)+2$. According to our integrability assumption, we can pick $q$ such that the $(k-1)q$-th moment of the steps of the random walk is finite. Note that $|\overline Y_{m,l}|\prod_{(i,j)\not=(m,m+1)}|X_i(l)-X_j(l)|l^{-\frac 12}$ can be estimated against a sum of products of the form 
$$
|Y_\gamma(l)|\prod_{i=1}^k \Big(\frac{|X_i(l)|}{\sqrt l}\Big)^{\beta_i},\qquad\mbox{for some }\gamma\in\{1,\dots,k\}\mbox{ and }\beta_1,\dots,\beta_k\in\N_0,
$$
with $\beta_1+\dots+\beta_k=k-1$. We use $C$ to denote a generic positive constant, not depending on $x$, nor on $l$, possibly changing its value from appearance to appearance. According to our moment assumption, for any $i\in\{1,\dots,k\}$, $l\in\N$ and any $\beta\in(0,q(k-1)]$, we have
\begin{equation} \label{momentnew}
\E_x\Big[\Big(\frac{|X_i(l)|}{\sqrt l}\Big)^{\beta}\Big]\leq C+\Big(\frac{|x|}{\sqrt l}\Big)^{\beta}\qquad\mbox{and}\qquad 
\E_x\Big[\Big(\frac{|X_i(l)|}{\sqrt l}\Big)^{\beta}|Y_i(l)|\Big]
\leq C+\Big(\frac{|x|}{\sqrt l}\Big)^{\beta}.
\end{equation}
On our set of summation on $l$, both upper bounds are bounded by $(C+ |x|^{-a})^\beta$ and may therefore be estimated against $C$, since $|x|$ is bounded away from zero for $x\in M$. Hence, the term on the last line of \eqref{Hoeldertrick} is bounded, since the sum over $l^{-q\xi[\frac k2(k-1)-1]}$ converges. 

The first term of the right hand side of \eqref{Hoeldertrick} can be estimated as follows.
$$
\begin{aligned}
\E_x&\Big[\tau^{p(\frac k2(k-1)-1)(\frac 12 +\xi)}\1_{\{|x|^{2(1+a)}<\tau\leq n\}}\Big]^{1/p}
= \Big( \int_{|x|^{2(1+a)(r+1)}}^{n^{r+1}} \d s \, \P_x(\tau^{r+1}>s) \Big)^{1/p} \\ &= C \Bigl( \int_{|x|^{2(1+a)}}^n \d t\, \P_x(\tau>t)t^r \Big)^{1/p}
\leq C \Big( \sum_{l={\lceil |x|^{2(1+a)} \rceil}}^n l^r \P_x(\tau >l) \Bigr)^{1/p}.
\end{aligned}
$$ 
Now put $\lambda = 1/p$. From this estimate, the assertion follows.
\qed
\end{proofsect}

We remark that upper estimates for $\E_x[|\Delta(X(\tau))|\1_{\{\tau=n\}}]$ in terms of $\P_x(\tau=n)$ are relatively easy to derive, but not sufficient for our purposes, since the techniques we develop in Section~\ref{sec-strategy} below do not give sufficient control on the asymptotics of $\P_x(\tau=n)$. However, they are good enough to control the ones of $\P_x(\tau>n)$.

\subsection{Local expansions.}\label{sec-expansions}
 
\noindent In this section we state, for future reference, one of our main tools, the expansion in the local central limit theorem, see \cite[Ch.~VII, Par. 3]{Pe75}. We are under the Centering and the Regularity Assumptions. Recall that, in the lattice case, the maximal span of the walk is $\alpha$, and in the non-lattice case that, for some $N\in\N$, the density $p_N$ is bounded. We define, for $n\geq N$,
\begin{equation}\label{pndef}
p_n(x)=\begin{cases}
\P_0(X_1(n)=x),&\mbox{in the lattice case,}\\
\frac{\P_0(X_1(n)\in\d x)}{\d x},&\mbox{in the non-lattice case.}
       \end{cases}
\end{equation}
In the lattice case, the numbers $p_n(x)$ sum up to one over the lattice $S=\alpha\Z$, and in the non-lattice case, $p_n$ is a probability density on $S=\R$. Then we have

\begin{lemma}[Local CLT expansion]\label{lem-LocCLT} Assume that the Centering and the Regularity Assumptions hold and that the $\mu$-th moment of the step distribution is finite for some $\mu \geq 2$. Then
\begin{equation}\label{lclt}
\sqrt n p_n(x\sqrt n)=\frac 1{\sqrt{2\pi }}{\rm e }^{-\frac 12 x^2}\big(1+o(1)\big)+\frac{o\big(n^{1-\frac\mu2}\big)}{1+|x|^\mu},\qquad\mbox{uniformly for $x\in \frac 1{\sqrt n}S$}.
\end{equation}
\end{lemma}

\begin{proofsect}{Proof} \cite[Thms.~VII.16 resp.~VII.17]{Pe75} state that 
\begin{equation}\label{locexp}
\sqrt n p_n(x\sqrt n)=\frac 1{\sqrt{2\pi }}{\rm e}^{-\frac 12 x^2}\Big(1+\sum_{\nu=1}^{\mu-2}\frac{\widetilde q_\nu(x)}{n^{\nu/2}}\Big)+\frac{o\big(n^{1-\frac\mu2}\big)}{1+|x|^\mu},\qquad\mbox{uniformly for $x\in \frac 1{\sqrt n}S$,}
\end{equation}
where $\widetilde q_\nu$ are polynomials of order $\leq 3\nu$ whose coefficients depend on the first cumulants of the step distribution only. The term ${\rm e}^{-\frac 12 x^2}\sum_{\nu=1}^{\mu-2}\frac{\widetilde q_\nu(x)}{n^{\nu/2}}$ is either equal to ${\rm e}^{-\frac 12 x^2}o(1)$ (if $|x|\leq o(n^{\frac 16})$) or it is $o(n^{1-\frac\mu2})$ (if $|x|\geq n^{\frac 17}$, say). Hence, \eqref{lclt} follows. 
\qed
\end{proofsect}

We are going to rephrase an integrated version of the Karlin-McGregor formula of Proposition~\ref{KMG} in terms of $p_n$. This will be the starting point of our proofs. For notational convenience, we will put $N=1 $ in the non-lattice case (it is easy, but notationally nasty, to adapt the proofs to other values of $N$). We introduce a rescaled version of the measure ${\Dcal}_{l}(x,\d y)$. For $x=(x_1,\dots,x_k)\in W\cap S^k$ and $l\in\N$, introduce
\begin{equation}\label{Dnscaldef}
\begin{aligned}
D_l^{\ssup n}(x,y)&=\det\Big[\Big(\sqrt n p_l\big(y_i\sqrt n-x_{j}\big)\Big)_{i,j=1,\dots,k}\Big]\\
&=\sum_{\sigma\in\Sym_k}\sign(\sigma) \prod_{i=1}^k\Big(\sqrt n p_l\big(y_i\sqrt n-x_{\sigma(i)}\big)\Big),\qquad y=(y_1,\dots,y_k)\in W.
\end{aligned}
\end{equation}
In the non-lattice case, the map $y\mapsto D_l^{\ssup n}(x,y)$ is a density of the measure that is obtained from ${\Dcal}_{l}(x,\d y)$ as the image measure under the map $y\mapsto y\sqrt n$. In the lattice case, it is equal to that measure; note that it is zero outside the lattice  $\frac 1{\sqrt n} S^k$.

Then we may rephrase Proposition~\ref{KMG} as follows.

\begin{lemma}\label{KMGrephrase} For any continuous and bounded function $f\colon \R^k\cap W\to\R$,
\begin{equation}\label{taurepr}
\begin{aligned}
\E_x&\Big[f\Big(n^{-\frac 12}X(n)\Big)\1_{\{\tau >n\}}\Big]\\
&=\begin{cases} \int_W \Big[D_n^{\ssup n}(x,y)-\E_x\big[\1_{\{\tau\leq n\}} D_{n-\tau}^{\ssup n}(X(\tau),y)\big]\Big]\, f(y)\,\d y&\mbox{in the non-lattice case,}\\
\sum_{y\in W\cap \frac 1{\sqrt n} S} \Big[D_n^{\ssup n}(x,y)-\E_x\big[\1_{\{\tau\leq n\}} D_{n-\tau}^{\ssup n}(X(\tau),y)\big]\Big]\, f(y)&\mbox{in the lattice case.}
\end{cases}
\end{aligned}
\end{equation}
\end{lemma}

\begin{proofsect}{Proof} This is a reformulation of \eqref{KMGneu} using \eqref{Dnscaldef}.
\qed
\end{proofsect}

\subsection{Our main strategy.}\label{sec-strategy}
 
\noindent We are going to explain how we will prove the integrability of $\Delta(X(\tau))$ under $\P_x$. As is seen from Lemma~\ref{integr}, our main task is to give good bounds for $\P_x(\tau>n)$, more precisely, to show that this quantity decays on the scale $n^{-\frac k4(k-1)}$ (where some error of some small positive power of $n$ would not spoil the proof). In order to do that, we use the Karlin-McGregor formula of Proposition~\ref{KMG}, more precisely, the formula in \eqref{taurepr}.

We shall need a cutting argument for large values of $y$ in \eqref{taurepr}.
To do this, we fix a slowly divergent scale function $n^\eta$ (with some small $\eta>0$) and cut the integral in 
\eqref{taurepr} into the area where $|y|\leq n^\eta$ and the remainder. The remainder is small, according to some moment estimate. On the set where $|y|\leq n^\eta$, we derive uniform convergence of the integrand, using appropriate expansions in the local central limit theorem and an expansion in the determinant. The second term in the integrand in \eqref{taurepr}, which is $\E_x\big[\1_{\{\tau\leq n\}} D_{n-\tau}^{\ssup n}(X(\tau),y)\big]$, will have to be split into the three parts where 
$$
\tau\leq t_n,\qquad t_n<\tau\leq n-s_n,\qquad n-s_n<\tau\leq n,
$$ 
where $t_n, s_n\to\infty$ are auxiliary sequences such that $n/t_n$ and $s_n/\sqrt n$ are small positive powers of $n$, say. The necessity for such a split is the application of the local central limit theorem inside the expectation: it is possible only for values $\tau\leq n-s_n$, and it gives asymptotically correct values only for $\tau\leq t_n$. We will eventually show that the first part gives the main contribution, and the two other parts are small remainder terms.

The reason that we have to work with a {\it local\/} central theorem (in fact, even with an expansion to sufficient deepness) is the following. After application of the approximation, there will be an extinction of terms in the determinant. As a result, the main term will turn out to be of order $n^{-\frac k4(k-1)}$. Hence, we need an error term of the size $o(n^{-\frac k4(k-1)})$ in the central limit regime, and this precision can be achieved only via an expansion of a local theorem. Sufficient moment conditions will imply that the contribution from outside the central limit region is also of the size $o(n^{-\frac k4(k-1)})$.

Let us now turn to the details. We first turn to the analysis of the main term of the integrand on the right of \eqref{taurepr}, where $\{\tau\leq n\}$ is replaced by $\{\tau\leq t_n\}$. We need the function 
\begin{equation}\label{Vndef}
V_n(x)=\Delta(x)-\E_x[\Delta(X(\tau))\1_{\{\tau\leq n\}}],\qquad n\in \N, \, x \in \R^k.
\end{equation}
Under the assumption that the steps have finite $(k-1)$-st moment (which we will in particular impose), it is clear that $V_n(x)$ is well-defined. It is also relatively easy to show that $V_n$ is positive:

\begin{lemma}\label{Vharm}  For any $n\in\N$ and $x\in W$, $V_n(x)>0$.
\end{lemma}

\begin{proofsect}{Proof} We recall the fact that $(\Delta(X(n)))_n$ is a martingale under $\P_x$ for any $x\in W$, see \cite[Th.~2.1]{KOR02}. Hence, we can calculate
$$
\begin{aligned}
V_n(x)&=\Delta(x)-\E_x\bigl[\1_{\{\tau\leq n\}}\E_{X(\tau)}[\Delta(X(n-m))]|_{\tau=m}\bigr]\\
&=\E_x[\Delta(X(n))]-\E_x[\Delta(X(n))\1_{\{\tau\leq n\}}]\\
&=\E_x[\Delta(X(n))\1_{\{\tau> n\}}],
\end{aligned}
$$
and this is obviously positive. 
\qed
\end{proofsect}

\begin{lemma}\label{CLTtool} Assume that the $\mu$-th moment of the steps is finite for some $\mu\geq k-1$. Fix small parameters $\eta,\xi_1>0$ satisfying $8\eta<\xi_1$. We put $t_n=n^{1-{\xi_1}}$. Then, for any $x,y\in W\cap\frac 1{\sqrt n}S^k$, uniformly for $|x|=o(\sqrt {t_n})$ and $|y|=o( n^\eta)$, as $n\to\infty$,
\begin{equation}\label{tau>ndistr}
\begin{aligned}
D_n^{\ssup n}(x,y)&-\E_x\big[\1_{\{\tau\leq t_n\}} D_{n-\tau}^{\ssup n}(X(\tau),y)\big]\\
&=\prod_{l=0}^{k-1}\frac 1{l!}n^{-\frac k4(k-1)} V_{t_n}(x)\, \frac{{\rm e}^{-\frac 12|y|^2}}{(2\pi)^{k/2}} \Delta(y)\,(1+o(1))+\Ocal(n^{1-\mu(\eta+{\xi_1}/4)}) .
\end{aligned}
\end{equation}
\end{lemma}

\begin{proofsect}{Proof} According to \eqref{lclt}, we have, uniformly in $x,y\in\R^k\cap\frac 1{\sqrt n}S^k$,
\begin{equation}\label{CLTDn}
\begin{aligned}
D_n^{\ssup n}(x,y)&=(2\pi)^{-\frac k2}\det\Big[\big({\rm e}^{-\frac 12 (y_j-x_i/\sqrt n)^2}\big)_{i,j=1,\dots,k}\Bigr] \,\big(1+o(1)\big)+o\big(n^{1-\frac\mu2}\big)\\
&=\frac{{\rm e}^{-\frac 12|y|^2}}{(2\pi)^{k/2}}{\rm  e}^{-\frac 1{2n}|x|^2}\det\Big[\big({\rm e}^{x_i y_j n^{-1/2}}\big)_{i,j=1,\dots,k}\Big] \,\big(1+o(1)\big)+o\big(n^{1-\frac\mu2}\big).
\end{aligned}
\end{equation}
In order to evaluate the last determinant, we write, for $|x|\,|y|=o(\sqrt n)$,
\begin{equation}\label{expseries}
{\rm e}^{x_i y_j n^{-1/2}}=\sum_{l=1}^k\frac {x_i^{l-1}}{(l-1)!}\,\frac{y_j^{l-1}}{\sqrt n^{l-1}}\Big[1+\Ocal\Big(\frac{|x|\,|y|}{\sqrt n}\Big)^k\Big],
\end{equation}
and use the determinant multiplication theorem. This gives, with $K'=\prod_{l=0}^{k-1}(l)!^{-1}$, if $|x|\,|y|=o(\sqrt n)$, 
\begin{equation}\label{deteval}
\begin{aligned}
\det\Bigl[({\rm e}^{x_i y_j n^{-1/2}})_{i,j=1,\dots,k}\Bigr]&=\det\Big[\Big(\frac {x_i^{l-1}}{(l-1)!}\Big)_{i,l=1,\dots,k}\Big]
\det\Big[\Big(\frac{y_j^{l-1}}{\sqrt n^{l-1}} \Big)_{l,j=1,\dots,k}\Big](1+o(1))\\
&=K'\Delta(x)n^{-\frac k4(k-1)}\Delta(y)(1+o(1)).
\end{aligned}
\end{equation}
Substituting \eqref{deteval} in \eqref{CLTDn}, we obtain
\begin{equation}\label{detident1a}
D_n^{\ssup n}(x,y)= K' n^{-\frac k4(k-1)}\Delta(x)\frac{{\rm e}^{-\frac 12|y|^2}}{(2\pi)^{k/2}}\Delta(y) \,(1+o(1))+o\big(n^{1-\frac \mu2}\big),\qquad \mbox{for }|x|\, |y|=o(\sqrt n).
\end{equation}

In order to handle the second term on the left of \eqref{tau>ndistr} we have to distinguish if $|X(\tau)|$ is large or not. For this purpose, fix $m_n=n^{1-{\xi_1}/2}$, and split
\begin{equation}\label{detident2}
\begin{aligned}
\E_x\big[&\1_{\{\tau\leq t_n\}} D_{n-\tau}^{\ssup n}(X(\tau),y)\big]\\
&=\E_x\Bigl[\1_{\{\tau\leq t_n\}}\1_{\{|X(\tau)|\leq n^\eta\sqrt {m_n}\}}D_{n-\tau}^{\ssup n}(X(\tau),y)\Bigr]+\E_x\Bigl[\1_{\{\tau\leq t_n\}}\1_{\{|X(\tau)|> n^\eta\sqrt {m_n}\}}D_{n-\tau}^{\ssup n}(X(\tau),y)\Bigr].
\end{aligned}
\end{equation}
Let us estimate the second term. From now we use $C$ to denote a generic positive constant, depending only on $k$ or the step distribution. 
Observe from \eqref{locexp} and \eqref{Dnscaldef} that
\begin{equation}\label{densbound}
\big|D_l^{\ssup n}(x,y)\big|\leq C\Big(\frac nl\Big)^{\frac k2},\qquad n,l\in\N,x,y\in W \cap \frac{1}{\sqrt{n}} S^k.
\end{equation}
Hence, on $\{\tau\leq t_n\}$, $D_{n-\tau}^{\ssup n}(X(\tau),y)$ is uniformly bounded in $n$ and $y$, since $t_n=o(n)$ and therefore $\frac n{n-\tau}$ is bounded. 
Using the boundedness of $D_{n-\tau}^{\ssup n}(X(\tau),y)$, the Markov inequality and \eqref{momentnew}, we obtain, for all $x,y\in W \cap \frac{1}{\sqrt{n}} S^k$ 
satisfying $|x|=o(\sqrt{t_n})$, as $n\to\infty$,
\begin{equation}\label{cutestiX(tau)}
\begin{aligned}
\Big|\E_x\Big[&\1_{\{\tau\leq t_n\}}\1_{\{|X(\tau)|> n^\eta\sqrt {m_n}\}}D_{n-\tau}^{\ssup n}(X(\tau),y)\Big]\Big|\\
&\leq C \sum_{l=1}^{t_n}\P_x\big(\tau=l,|X(l)|> n^\eta\sqrt {m_n}\big)
\leq  C t_n \sup_{l=1}^{t_n}\P_x(|X(l)|> n^\eta\sqrt {m_n})\\
&\leq C t_n \sup_{l=1}^{t_n}\E_x[|X(l)|^\mu]n^{-\eta\mu} m_n^{-\mu/2}
\leq \Ocal\big(t_n (t_n/m_n)^{\mu/2}n^{-\eta\mu}\big)\\
&\leq \Ocal\big(n^{1-\mu(\eta+{\xi_1}/4)}\big).
\end{aligned}
\end{equation}
The first term in \eqref{detident2} can be handled in the same way as in \eqref{detident1a}. Indeed, for $x,y\in W \cap \frac{1}{\sqrt{n}} S^k$ satisfying $|x|=o(\sqrt{t_n})$ and $|y|=o(n^\eta)$,
\begin{equation}\label{detident3}
\begin{aligned}
\E_x\Bigl[&\1_{\{\tau\leq t_n\}}\1_{\{|X(\tau)|\leq n^\eta\sqrt {m_n}\}}D_{n-\tau}^{\ssup n}(X(\tau),y)\Bigr]\\
&= \prod_{l=0}^{k-1}\frac1{l!}n^{-\frac k4(k-1)}\E_x\big[\1_{\{\tau\leq t_n\}}{\rm e}^{-\frac 1{2n}|X(\tau)|^2}\1_{\{|X(\tau)|\leq n^\eta\sqrt {m_n}\}\}}\Delta(X(\tau))\big]\frac{{\rm e}^{-\frac 12|y|^2}}{(2\pi)^{k/2}}\Delta(y) \,(1+o(1))\\
&\qquad\qquad\qquad\qquad\qquad+o\big(n^{1-\frac\mu2}\big)\\
&=\prod_{l=0}^{k-1}\frac 1{l!} n^{-\frac k4(k-1)}\E_x\big[\1_{\{\tau\leq t_n\}}\Delta(X(\tau))\big]\frac{{\rm e}^{-\frac 12|y|^2}}{(2\pi)^{k/2}}\Delta(y) \,(1+o(1))+o\big(n^{1-\frac\mu2}\big).
\end{aligned}
\end{equation}
The first step uses that \eqref{detident1a} is applicable for $x=X(\tau)$ since $|X(\tau)|\,|y|=o(n^{2\eta}\sqrt {m_n})=o(n^{\frac 12-\frac{\xi_1}4+2\eta})=o(\sqrt n)$; recall that we assumed that $8\eta <\xi_1$. The second step in \eqref{detident3} is derived in a similar way as in \eqref{cutestiX(tau)}, using also that, on the event $\{ |X(\tau)|\leq n^\eta\sqrt {m_n}\}$, we have $|X(\tau)|^2/n\leq n^{2\eta} m_n/n\to 0$, according to the choice of $m_n$ and $8 \eta < \xi_1$.

Now substitute \eqref{detident3} and \eqref{cutestiX(tau)} in \eqref{detident2} and this and \eqref{detident1a} on the left side of \eqref{tau>ndistr} 
to finish the proof.
\end{proofsect}
\qed

Now we examine the part where $t_n\leq \tau\leq n-s_n$. We restrict to the non-lattice case; the formulation for the lattice case and its proof are analogous, and are left to the reader.

\begin{lemma}\label{middleterm} Consider the non-lattice case. Assume that the $\mu$-th moment of the steps is finite for some $\mu\geq k-1$. Fix small 
parameters $\eta,\eps,{\xi_1},\xi_2>0$ such that $\xi_2>\eps+\eta$. We put $t_n=n^{1-{\xi_1}}$ and $s_n=n^{\frac 12+\xi_2}$. 
Then, for any $x\in W$, uniformly for $|x|=o(\sqrt n)$, as $n\to\infty$,
\begin{equation}\label{esti2term}
\begin{aligned}
\int_W&\1_{\{|y|\leq n^\eta\}}\Big|\E_x\big[\1_{\{t_n\leq \tau\leq n-s_n\}} D_{n-\tau}^{\ssup n}(X(\tau),y)\big]\Big|\,\d y\\
&\leq  \Ocal\big(n^{\eps+\eta-\xi_2}\big)\P_x(\tau\geq t_n)+o\big(n^{\frac 12-\frac \mu 4+\eta k+\xi_2}\big)+\Ocal \big(n^{1-\mu(\eta+{\xi_1}/4)}\big)+\Ocal\big(n^{1-\eps\mu}\big).
\end{aligned}
\end{equation}
\end{lemma}

\begin{proofsect}{Proof} For notational convenience, we assume that $t_n$ and $s_n$ are integers. We will use $C>0$ as a generic positive constant, depending only on the step distribution or on $k$, possibly changing its value from appearance to appearance.

Similarly to \eqref{detident2}, in the expectation on the left of \eqref{esti2term}, we distinguish if $|X(\tau)|$ is larger than $n^\eta \sqrt {m_n}$ or not, where this time we pick $m_n=n^{1+{\xi_1}/2}$. Furthermore, we recall the stopping time $\tau_{m}$ from \eqref{tauijdef} and sum on all values of $m$ and distinguish if $|X_{m}(\tau)-X_{m+1}(\tau)|$ is smaller than $n^\eps$ or not. Furthermore, we sum on all values of $\tau$. 
Recalling also \eqref{densbound}, we estimate the integrand on the left side of \eqref{esti2term} as
\begin{equation}\label{estimiddlesplit}
\Big|\E_x\big[\1_{\{t_n\leq \tau\leq n-s_n\}} D_{n-\tau}^{\ssup n}(X(\tau),y)\big]\Big|
\leq Z_n^{\ssup 1}(y)+Z_n^{\ssup 2}(y)+Z_n^{\ssup 3}(y),
\end{equation}
where 
\begin{eqnarray}
Z_n^{\ssup 1}(y)&=&\sum_{l=t_n}^{n-s_n}\E_x\Big[\1_{\{\tau=l\}}\1_{\{|X(l)|>n^\eta\sqrt {m_n}\}}\big|D_{n-l}^{\ssup n}(X(l),y)\big|\Big],\\
Z_n^{\ssup 2}(y)&=&\sum_{l=t_n}^{n-s_n}\sum_{m=1}^{k-1} \E_x\Big[\1_{\{\tau=\tau_{m}=l\}}\1_{\{|X_{m}(l)-X_{m+1}(l)|>n^\eps\}}\big|D_{n-l}^{\ssup n}(X(l),y)\big|\Big]\\
Z_n^{\ssup 3}(y)&=&\sum_{l=t_n}^{n-s_n}\sum_{m=1}^{k-1} \E_x\Big[\1_{\{\tau=\tau_{m}=l\}}\1_{\{|X(l)|\leq n^\eta\sqrt {m_n}\}}\1_{\{|X_{m}(l)-X_{m+1}(l)|\leq n^\eps\}} \big|D_{n-l}^{\ssup n}(X(l),y)\big|\Big].\label{Z3def}
\end{eqnarray}

Let us estimate the integral over $Z_n^{\ssup 1}$. Observe that the measure with density $y\mapsto |D_{n-l}^{\ssup n}(x,y)|$ has bounded total mass:
\begin{equation}\label{bdtotalmass}
\int_W  \d y\,\big|D_{n-l}^{\ssup n}(x,y)\big|\leq k! \Big(\int_{\R}\d z\, \sqrt n p_{n-l}\big(z\sqrt n\big)\Big)^k\leq k!,\qquad 0\leq l\leq n, x\in W.
\end{equation}
Using this and  the Markov inequality, we obtain, for $|x|=o(\sqrt n)=o(\sqrt{n-s_n})$, as $n\to\infty$,
\begin{equation}\label{Z1esti}
\begin{aligned}
\int_W\d y\,Z_n^{\ssup 1}(y)&\leq C\sum_{l=t_n}^{n-s_n}\P_x\big(|X(l)|>n^\eta\sqrt {m_n}\big)
\leq C \sum_{l=t_n}^{n-s_n}\E_x\big[|X(l)|^\mu\big]n^{-\eta\mu} m_n^{-\mu/2}\\
&\leq C n (n/m_n)^{\mu/2}n^{-\eta\mu}\leq \Ocal\big(n^{1-\mu(\eta+{\xi_1}/4)}\big).
\end{aligned}
\end{equation}

In a similar way, we estimate the integral over $Z_n^{\ssup 2}$. Recall that, on the event $\{\tau=l=\tau_{m}\}$, we may use \eqref{overshoot}, \eqref{bdtotalmass} and then the Markov inequality, to obtain
\begin{equation}\label{Z2esti}
\int_W\d y\,Z_n^{\ssup 2}(y)\leq C \sum_{l=t_n}^{n-s_n}\sum_{m=1}^{k-1}\P_x\big(|Y_{m+1}(l)-Y_{m}(l)|>n^\eps\big)
\leq Cn  n^{-\mu \eps}\leq \Ocal\big(n^{1-\mu\eps}\big).
\end{equation}
Hence, \eqref{Z1esti} and \eqref{Z2esti} yield the two latter error terms on the right hand side of \eqref{esti2term}.

We turn to an estimate of $Z_n^{\ssup 3}(y)$. For $t_n\leq l\leq n-s_n$, on the event $\{\tau=l\}\cap\{|X(l)|\leq n^\eta \sqrt {m_n}\}$, we use the local central limit theorem in \eqref{lclt} to obtain
\begin{equation}\label{estitntaun-sn1}
\begin{aligned}
D_{n-l}^{\ssup n}&(X(l),y)\\
&=\det\Big[\Big(\sqrt n p_{n-l}\big(y_j\sqrt n-X_i(l)\big)\Big)_{i,j=1,\dots,k}\Big]\\
&=\Big(\frac n{n-l}\Big)^{\frac k2}
\frac{{\rm e}^{-\frac 12|y|^2\frac n{n-l}}}{(2\pi)^{k/2}} {\rm e}^{-\frac 12\frac 1{n-l} |X(l)|^2}
\det\Big[\Big({\rm e}^{y_j X_i(l)\frac{\sqrt n}{n-l}}\Big)_{i,j=1,\dots,k}\Big] \,(1+o(1))+o\big((n-l)^{1-\mu/2}\big).
\end{aligned}
\end{equation}
Abbreviate $a=X(l)\sqrt n/(n-l)$, then the determinant may be written
\begin{equation}\label{detrewrite}
\det\Big[\Big({\rm e}^{a_i y_j}\Big)_{i,j=1,\dots,k}\Big]=\sum_{1\leq i<j\leq k}\Big(1-{\rm e}^{(a_m-a_{m+1})(y_i-y_j)}\Big)\sum_{\heap{\sigma\in\Sym_k}{\sigma(m)=i,\sigma(m+1)=j}}\sign(\sigma)\prod_{i=1}^k{\rm e}^{a_{\sigma(i)}y_i}.
\end{equation}
Observe that the exponent $(a_m-a_{m+1})(y_i-y_j)$ asymptotically vanishes on the event $\{|X_{m}(l)-X_{m+1}(l)|\leq n^\eps\}$, since on $\{|y| \leq n^{\eta} \}$ 
we obtain
\begin{equation}\label{asmall}
|a_m-a_{m+1}|\,|y_i-y_j|\leq n^{\eps+\eta}\frac{{\sqrt n}}{n-l}\leq n^{\eps+\eta+\frac 12}\frac 1{s_n}=n^{\eps+\eta-\xi_2}=o(1),
\end{equation}
since $\xi_2>\eps+\eta$. 
Hence, we may use that $|1-{\rm e}^x|\leq \Ocal(|x|)$ as $x\to0$, and obtain from \eqref{detrewrite} that 
$$
\begin{aligned}
\Big|\det\Big[\Big({\rm e}^{y_j X_i(l)\frac{\sqrt n}{n-l}}\Big)_{i,j=1,\dots,k}\Big]\Big|
&\leq \Ocal\big(n^{\eta+\eps -\xi_2}\big)\sum_{\sigma\in\Sym_k}\prod_{i=1}^k{\rm e}^{X_{\sigma(i)}(l)\frac{\sqrt n}{n-l}y_i}.
\end{aligned}
$$
Using this in \eqref{estitntaun-sn1} and this in \eqref{Z3def}, we arrive at 
$$
Z_n^{\ssup 3}(y)\leq \Ocal\big(n^{\eta+\eps -\xi_2}\big)\sum_{l=t_n}^{n-s_n}\Big(\frac n{n-l}\Big)^{\frac k2}\sum_{\sigma\in\Sym_k}\E_x\Big[\1_{\{\tau=l\}}{\rm e}^{-\frac 12|X_\sigma(l)/\sqrt n-y|^2\frac n{n-l}} \Big]+o\big(s_n^{1-\mu/2}\big).
$$
Now we integrate over $y$ and use Fubini's theorem:
\begin{equation}\label{Z3esti2}
\begin{aligned}
\int_W &\1_{\{|y|\leq n^\eta\}}Z_n^{\ssup 3}(y)\,\d y\\
&\leq \Ocal\big(n^{\eta+\eps -\xi_2}\sum_{l=t_n}^{n-s_n} \sum_{\sigma\in\Sym_k}\E_x\Big[\1_{\{\tau=l\}}\int_{\R^k}{\rm e}^{-\frac 12|X_\sigma(l)/\sqrt n-y|^2\frac n{n-l}}\Big(\frac n{n-l}\Big)^{\frac k2}\,\d y \Big]+o\big(n^{\eta k}s_n^{1-\mu/2}\big)\\
&\leq \Ocal\big(n^{\eps+\eta-\xi_2}\big)\sum_{l=t_n}^{n-s_n} \P_x(\tau=l)+o\big(n^{\frac 12-\frac \mu 4+\eta k+\xi_2}\big)\\
&\leq \Ocal\big(n^{\eps+\eta-\xi_2}\big)\P_x(\tau\geq t_n)+o\big(n^{\frac 12-\frac \mu 4+\eta k+\xi_2}\big).
\end{aligned}
\end{equation}
Substituting this in \eqref{estimiddlesplit} and combining with \eqref{Z1esti} and \eqref{Z2esti}, we arrive at \eqref{esti2term}.
\qed
\end{proofsect}

\subsection{Proof of integrability of $\boldsymbol{\Delta(X(\tau))}$.}\label{sec-integrability}

\noindent Now we collect the preceding and prove the integrability of $\Delta(X(\tau))$ under $\P_x$ for any $x$:

\begin{prop}[Integrability of ${\Delta(X(\tau))}$]\label{Xtaukint} There is $\mu_k$, depending only on $k$, such that, if the $\mu_k$-th moment of the steps is finite, for any $x\in W$, the variable $\Delta(X(\tau))$ is integrable under $\P_x$. Moreover, uniformly in $x\in W$ on compacts,
\begin{equation}\label{tauX(n)asy}
\lim_{n\to\infty}n^{\frac k4(k-1)}\P_x\Big(\tau>n;\frac 1{\sqrt n}X(n)\in\d y\Big)=\prod_{l=1}^{k-1}\frac 1{l!}\,V(x)\frac{{\rm e}^{-\frac 12 |y|^2}}{(2\pi)^{k/2}}\Delta(y)\,\d y,
\end{equation}
in the sense of weak convergence of distributions on $W$.
\end{prop}
 
Before we give the proof, we state some conclusions:

\begin{cor}[Tails of $\tau$ and limit theorem before violation of ordering]\label{cor-tautail} Assume that the $\mu_k$-th moment of the steps is finite, with $\mu_k$ as in Proposition~\ref{Xtaukint}. Then
\begin{equation}\label{tauasy}
\lim_{n\to\infty}n^{\frac k4(k-1)}\P_x(\tau>n)=KV(x), \qquad\mbox{where } K=\prod_{l=1}^{k-1}\frac 1{l!}\,\int_W\frac{{\rm e}^{-\frac 12 |y|^2}}{(2\pi)^{k/2}}\Delta(y)\,\d y.
\end{equation}
Furthermore, the distribution of $n^{-\frac 12}X(n)$  under $\P_x(\,\cdot\mid\tau>n)$ converges towards the distribution on $W$ with density $y\mapsto\frac 1{Z_1}{\rm e}^{-\frac 12 |y|^2}\Delta(y)$, where $Z_1$ is the normalisation.
\end{cor}

Note that we {\it a priori\/} do not know yet whether or not $V$ is positive. However, its nonnegativity is now clear, which we want to state explicitly:

\begin{cor}[Nonnegativity of $V$]\label{Vpos} For any $x\in W$, the number $V(x)$ defined in \eqref{Vdef} is well-defined and nonnegative.
\end{cor}

This follows either from the asymptotics in Proposition~\ref{Xtaukint} or from the fact that $V_n$ is positive (see Lemma~\ref{Vharm}), in combination with $V(x)=\lim_{n\to\infty}V_n(x)$ via Lebesgue's theorem.

\begin{proofsect}{Proof of Proposition~\ref{Xtaukint}} Fix some continuous and bounded function $f\colon W\to\R$. We abbreviate 
$$
a_n(f)=n^{\frac k4(k-1)}\E_x\Big[f\Big(n^{-\frac 12}X(n)\Big)\1_{\{\tau>n\}}\Big]\qquad \mbox{and}\qquad a_n=a_n(\1) \qquad\mbox{and}\qquad A_n=\max\{a_1,\dots,a_n\}.
$$
Our first step, see \eqref{anesti2}, is to derive an expansion for $a_n(f)$ in terms of $K(f) V_{t_n}(x)$ with some suitable $K(f)\in(0,\infty)$ and error terms depending on $A_n$ and $\P_x(n-s_n\leq \tau\leq n)$. Specialising to $f=\1$ and using an upper bound for $V_{t_n}(x)$ derived from Lemma~\ref{integr}, we obtain a recursive upper bound for $a_n$ in terms of $A_n$ and $a_{n-s_n}-a_n$. This estimate directly implies that $(A_n)_{n\in\N}$ is bounded, hence also $(a_n)_{n\in\N}$ is bounded. Via Lemma~\ref{integr}, this directly implies the integrability of $\Delta(X(\tau))$, i.e., we know that $V(x)$ is well-defined, and $V(x)=\lim_{n\to\infty}V_{t_n}(x)$. Using this again in \eqref{anesti2} for $f=\1$, we further obtain that $a_n$ converges towards $K V(x)$, where $K$ is defined in \eqref{tauasy}. Using this in \eqref{anesti2} for arbitrary $f$, we derive the assertions in \eqref{tauX(n)asy} and finish the proof.

As in Lemmas~\ref{CLTtool} and \ref{middleterm}, we pick small parameters $\eps,\eta,{\xi_1},\xi_2>0$, and we put $t_n=n^{1-{\xi_1}}$ and $s_n=n^{\frac 12+\xi_2}$. Now we require that
\begin{equation}\label{constantschoice}
8\eta<\xi_1\qquad\mbox{and}\qquad\xi_2>\eps+\eta+\xi_1\frac k4(k-1),
\end{equation}
and we pick $\mu$ so large that 
\begin{equation}\label{muchoice}
\mu > (k-1) \Big( \frac k2 (k-1) +2 \Big) \,\, \text{and} \,\, 
\frac k4(k-1) +\max\Big\{-\mu\eta,\frac 12-\frac\mu4 +\eta k+\xi_2, 1-\mu(\eta+\xi_1/4),1-\eps\mu \Big\}<0.
\end{equation}
In the following, we will restrict ourselves to the non-lattice case. The necessary changes for the lattice cases are only notational. We use $C$ to denote a generic positive constant that is uniform in $n$ and uniform in $x$ on compacts and may change its value from appearance to appearance. 
All following limiting assertions hold for $x\in W$  uniformly on compacts. We begin with multiplying \eqref{taurepr} with $n^{\frac k4(k-1)}$ and distinguishing the events $\{|X_n|\geq n^{\frac 12 +\eta}\}$ and its complement. This gives
\begin{equation}\label{ansplit}
\begin{aligned}
a_n(f)&= n^{\frac k4(k-1)}\E_x\Big[f\Big(n^{-\frac 12}X(n)\Big)\1_{\{\tau>n\}}\1_{\{n^{-1/2}|X(n)|\geq n^\eta\}} \Big]\\
&\qquad +n^{\frac k4(k-1)}\int_W\1_{\{|y|\leq n^\eta\}} \big[I_n(y)+II_n( y)+III_n(y)\big]\, f(y)\,\d y,
\end{aligned}
\end{equation}
where 
\begin{eqnarray}
I_n(y)&=&D_n^{\ssup n}(x,y)-\E_x\big[\1_{\{\tau\leq t_n\}} D_{n-\tau}^{\ssup n}(X(\tau),y)\big],\\
II_n(y)&=&-\E_x\big[\1_{\{t_n< \tau\leq n-s_n\}} D_{n-\tau}^{\ssup n}(X(\tau),y)\big],\\
III_n(y)&=&-\E_x\big[\1_{\{n-s_n< \tau\leq n\}} D_{n-\tau}^{\ssup n}(X(\tau),y)\big].
\end{eqnarray}

We use the Markov inequality for the first term on the right hand side of \eqref{ansplit}, which gives
\begin{equation}\label{integr1}
\begin{aligned}
n^{\frac k4(k-1)}&\Big|\E_x\Big[f\Big(n^{-\frac 12}X(n)\Big)\1_{\{\tau>n\}}\1_{\{n^{-1/2}|X(n)|\geq n^\eta\}} \Big]\Big|\\
&\leq n^{\frac k4(k-1)}\|f\|_\infty \P_x(
n^{-1/2}|X(n)|\geq n^\eta)
\leq C n^{\frac k4(k-1)}  n^{-\eta\mu}\\
&\leq \Ocal\big(n^{\frac k4(k-1)-\mu\eta}\big)=o(1),
\end{aligned}
\end{equation}
since $\frac k4(k-1)-\mu\eta<0$ by \eqref{muchoice}.
Hence, 
\begin{equation}\label{ansplit1}
a_n(f)=o(1)+\int_W\1_{\{|y|\leq n^\eta\}} n^{\frac k4(k-1)}\big[I_n(y)+II_n(y)+III_n(y)\big]\,f(y)\,\d y.
\end{equation}

Now we use \eqref{tau>ndistr} for $I_n(y)$, \eqref{esti2term} for $II_n(y)$ and \eqref{bdtotalmass} for $III_n( y)$. This gives
\begin{equation}\label{I-II-IIIesti}
\begin{aligned}
a_n(f)&=o(1)+\prod_{l=0}^{k-1}\frac 1{l!}\int_W\1_{\{|y|\leq n^\eta\}}\frac{{\rm e}^{-\frac 12 |y|^2}}{(2\pi)^{k/2}}\Delta(y)f(y) \,\d y\,V_{t_n}(x)(1+o(1))\\
&\qquad+ n^{\frac k4(k-1)}\Ocal\big(n^{\eps+\eta-\xi_2}\big)\P_x(\tau\geq t_n)+n^{\frac k4(k-1)}\Ocal\big(\P_x(n-s_n\leq \tau\leq n)\big) \\
&\qquad +n^{\frac k4(k-1)}\Big(o\big(n^{\frac 12-\frac\mu 4+\eta k+\xi_2}\big)+\Ocal(n^{1-\mu(\eta+{\xi_1}/4)}) +\Ocal\big(n^{1-\eps\mu}\big)\Big).
\end{aligned}
\end{equation}
Use that $ \lim_{n\to\infty}(n/t_n)^{\frac k4(k-1)}n^{\eps+\eta-\xi_2}= 0$ (see \eqref{constantschoice}) to write the first term in the second line of \eqref{I-II-IIIesti} as $o(1) a_{t_n}$. Using our assumptions on $\mu$ in \eqref{muchoice}, the third line of \eqref{I-II-IIIesti} is $o(1)$.  Therefore, we have from \eqref{I-II-IIIesti}
\begin{equation}\label{anesti2}
a_n(f) =K(f)\,V_{t_n}(x)(1+o(1))+o(1)+o(1)a_{t_n}+n^{\frac k4(k-1)}\Ocal\big(\P_x(n-s_n\leq \tau\leq n)\big),
\end{equation}
where 
$$
K(f)=\prod_{l=1}^{k-1}\frac 1{l!}\,\int_W\d y\,\frac{{\rm e}^{-\frac 12 |y|^2}}{(2\pi)^{k/2}}\Delta(y)f(y).
$$

In order to estimate $V_{t_n}(x)$, we use Lemma~\ref{integr} (more precisely, \eqref{summable}) and obtain, 
for some $\lambda\in(0,1)$ and $r\in (0,\frac k4(k-1)-1)$,
\begin{equation}
|V_{t_n}(x)|\leq |\Delta(x)|+C+C\Big(\sum_{l=1}^{t_n}l^r\P_x(\tau>l)\Big)^\lambda
\leq C+C\Big(\sum_{l=1}^{t_n}l^{r-\frac k4(k-1)}a_l\Big)^\lambda
\leq C+C A_n^\lambda,
\end{equation}
since the sum over $l^{r-\frac k4(k-1)}$ converges. Remark, that $\mu$ is chosen so large that \eqref{muchoice} holds, 
therefore Lemma~\ref{integr} can be applied.
We therefore obtain from \eqref{anesti2}, specialised to $f=\1$,
\begin{equation}
a_n\leq C +C A_n^\lambda+ o(1) A_n+C\big(a_{n-s_n}(1+o(1))-a_n\big),
\end{equation}
where we used $n^{k/4(k-1)} \P_x(n - s_n \leq \tau \leq n) = a_{n -s_n} (n/(n-s_n))^{\frac k4 (k-1)} - a_n$
and where we recalled that $s_n=o(n)$ and therefore $n^{\frac k4(k-1)}=(n-s_n)^{\frac k4(k-1)}(1+o(1))$.

Solving this recursive bound for $a_n$, we obtain the estimate
\begin{equation}\label{anesti1}
a_n\leq  C +\frac C{C+1} A_n^\lambda+\frac {C+o(1)}{C+1}a_{n-s_n}\leq  C +\frac C{C+1} A_n^\lambda+\frac {C+o(1)}{C+1}A_{n}.
\end{equation}
Since the right hand side is increasing in $n$, we also have this upper bound for $A_n$ instead of $a_n$. Now it is clear that $(A_n)_{n\in\N}$ is a bounded sequence, i.e., $\limsup_{n\to\infty}n^{\frac k4(k-1)}\P_x(\tau>n)<\infty$.
From Lemma~\ref{integr} (more precisely, from \eqref{summable}) we see that $V(x)$ is well-defined, and $\lim_{n\to\infty}V_n(x)=V(x)$ by Lebesgue's theorem.

Now, we prove that $\lim_{n\to\infty}a_n= K V(x)$, where $K=K(\1)$ is defined in \eqref{tauasy}. We start from \eqref{anesti2} with $f=\1$, which now reads
\begin{equation}\label{anrepr}
a_n= KV(x)\,(1+o(1))+o(1)+\Ocal\big(a_{n-s_n}(1+o(1))-a_n\big),\qquad n\to\infty,
\end{equation}
where we recall that $a_{n-s_n}(1+o(1))-a_n\geq 0$. Hence, 
$$
a_n\leq \frac{K V(x)}{1+C}(1+o(1))+\frac{C+o(1)}{1+C} a_{n-s_n}.
$$
Consider $\overline a=\limsup_{n\to\infty}a_n$ and note that $\limsup_{n \to \infty} a_{n -s_n} \leq \overline a$.
Hence, we obtain that $\overline a\leq KV(x)$. In the same way, we derive from \eqref{anrepr} that $\underline a=\liminf_{n\to\infty}a_n$ 
satisfies $\underline a\geq KV(x)$. This shows that $\lim_{n\to\infty}a_n$ exists and is equal to $K V(x)$. 

Now we go back to \eqref{anesti2} with an arbitrary continuous and bounded  function $f$, which now reads
$$
\begin{aligned}
a_n(f)& =K(f)\,V(x)(1+o(1))+o(1)+n^{\frac k4(k-1)}\Ocal\big(\P_x(n-s_n\leq \tau\leq n)\big)\\
&=K(f)\,V(x)(1+o(1))+o(1)+\Ocal\big(|a_{n-s_n}-a_n|\big)\\
&=K(f)\,V(x)(1+o(1))+o(1).
\end{aligned}
$$
This implies that \eqref{tauX(n)asy} holds and finishes the proof of Proposition~\ref{Xtaukint}.
\qed
\end{proofsect}

\section{Properties of $V$, and convergence towards Dyson's Brownian motions}\label{sec-Vpos}

\noindent In this section we prove a couple of properties of the function $V$ established in Proposition~\ref{Xtaukint}, like integrability properties, regularity, positivity and behavior at infinity. Furthermore, we introduce the Doob $h$-transform with $h=V$ and show that it is equal to a family of $k$ independent random walks, conditioned to be ordered at any time, and we prove an invariance principle towards Dyson's Brownian motions.
 
In Section~\ref{sec-Vexist} we saw that establishing the existence of $V(x)$ requires a detailed analysis of the asymptotics of $\P_x(\tau>n;n^{-\frac 12}X(n)\in\d y)$ on the scale $n^{-\frac k4(k-1)}$. In the present section, we will see that the proof of many properties of $V$, in particular its asymptotics at infinity and its positivity, requires a good control of the asymptotics of $\P_{x\sqrt n}(\tau>l;l^{-\frac 12}X(l)\in\d y)$ for $l/n$ bounded away from zero. 

Recall the notation of Section~\ref{sec-DysonBM}, in particular that a tuple of $k$ Brownian motions starts in $x\in\R^k$ under ${\tt P}_x$, and their collision time, $T$, in \eqref{Tdef}. 

\begin{lemma}\label{lem-Invariance} For any $t>0$ and any $x\in W$, 
\begin{equation}\label{DonskerDyson}
\lim_{n\to\infty}\P_{x\sqrt n}(\tau>tn;n^{-\frac 12}X(\lfloor tn\rfloor )\in\d y)={\tt P}_x(T>t; B(t)\in\d y)\qquad\mbox{weakly.}
\end{equation}
\end{lemma}
\begin{proofsect}{Proof} This follows from Donsker's invariance principle. Indeed, the vector $B^{\ssup{n}}=(B_1^{\ssup{n}},\dots,B_k^{\ssup{n}})$ of 
the processes $B_i^{\ssup{n}}(t)=n^{-\frac 12} X_i(\lfloor tn\rfloor)$ converges towards a Brownian motion $B$ on $\R^k$ starting at $x$. Since the event
$\{\forall s\in[0,t]\colon B^{\ssup{n}}_s\in W\}$ is open and $\{\forall s\in[0,t]\colon B^{\ssup{n}}_s\in \overline{W}\}$ is closed, we obtain
$\liminf_{n \to \infty} \P_x \big( \forall s\in[0,t]\colon B^{\ssup{n}}_s\in W \big) \geq \P_x \big( \forall s\in[0,t]\colon B_s\in W\big)$
and $\limsup_{n \to \infty} \P_x \big( \forall s\in[0,t]\colon B^{\ssup{n}}_s\in W \big) \leq \P_x \big( \forall s\in[0,t]\colon B_s\in \overline{W} \bigr)$.
Using that $\P_x \big( \forall s\in[0,t]\colon B_s\in W \big) = \P_x \big( \forall s\in[0,t]\colon B_s\in \overline{W} \big)$,  
the event $\{\tau>tn\}=\{\forall s\in[0,t]\colon B^{\ssup{n}}_s\in W\}$ converges in distribution to the event $\{\forall s\in[0,t]\colon B_s\in W\}=\{T>t\}$. Hence, we obtain \eqref{DonskerDyson}. 
\qed
\end{proofsect}

The following estimate is analogous to Lemma~\ref{integr}.

\begin{lemma}\label{integrneu}Assume that the $\mu$-th moment of the steps is finite, for some $\mu>(k-1)(\frac k2(k-1)+2)$. Then there are $\eps>0$, $r\in(0,\frac k4(k-1)-1)$ and $\lambda\in(0,1)$ such that, for any $M>0$, there is $C_M>0$ such that, for any $x\in W$ satisfying $ |x|\leq M$,
$$
\E_{\sqrt n\, x}[|\Delta(X(\tau))|\1_{\{\tau> Rn\}}]\leq C_M (Rn)^{-\eps}\Big(\sum_{l\in\N}l^r \P_{\sqrt n\, x}(\tau>l)\Big)^\lambda,\qquad R>1, n\in\N.
$$
\end{lemma}

\begin{proofsect}{Proof} Recall the notation of $\tau_m$ in \eqref{tauijdef} and the estimate in \eqref{overshoot}. In the same way as in the proof of Lemma~\ref{integr} (see in particular \eqref{Hoeldertrick}), we obtain the estimate, for any $m\in\{1,\dots,k-1\}$,
\begin{equation}\label{Hoeldertrickneu}
\begin{aligned}
\E_{\sqrt n\, x}&\big[|\Delta(X(\tau))|\1_{\{\tau=\tau_m> Rn\}}\big]\leq 
\E_{\sqrt n\, x}\big[\tau^{p[\frac k2(k-1)-1](\frac 12+\xi)}\1_{\{\tau\geq Rn\}}\big]^{1/p}\\
&\qquad\times 
\Big(\sum_{l\geq Rn}l^{-q \xi[\frac k2(k-1)-1]} \E_{\sqrt n\, x}\Big[|\overline Y_{m,l}|^q \prod_{(i,j)\not=(m,m+1)}\Big|\frac{X_i(l)-X_j(l)}{\sqrt l}\Big|^q\Big]\Big)^{1/q},
\end{aligned}
\end{equation}
where we pick again 
$\xi=(\frac k2(k-1)-1)^{-1}(\frac k2(k-1)+2)^{-1}$ and any $q>(\xi[\frac k2(k-1)-1])^{-1}=\frac k2(k-1)+2$; and $p$ is determined by $1=\frac 1p+\frac 1q$.

In a similar way as in the proof of Lemma~\ref{integr}, one sees that the expectation in the second line of \eqref{Hoeldertrickneu} is bounded in $n$ and $x\in W$ satisfying $|x|\leq M$, since $l/n$ is bounded away from zero in the sum. Hence, the second line is not larger than 
$$
C_M\Big(\sum_{l\geq Rn}l^{-q \xi[\frac k2(k-1)-1]}\Big)^{1/q}
\leq C_M (Rn)^{\frac 1q-\xi[\frac k2(k-1)-1]},
$$
for some $C_M>0$, not depending on $n$ nor on $R$ nor on $x$, as long as $x\in W$ and $|x|\leq M$. Note that the exponent $-\eps=\frac 1q-\xi[\frac k2(k-1)-1]$ is negative.

Let us turn to the other term on the right side of \eqref{Hoeldertrickneu}. We abbreviate $r=p[\frac k2(k-1)-1](\frac 12+\xi)-1$ and know that $r+1<\frac k4(k-1)$. Then we have
\begin{equation} \label{schoen}
\E_{\sqrt n\, x}\big[\tau^{p[\frac k2(k-1)-1](\frac 12+\xi)}\1_{\{\tau\geq Rn\}}\big]^{1/p}
\leq \E_{\sqrt n\, x}\big[\tau^{r+1}\big]^{1/p}\leq C\Big(\sum_{l\in\N}l^r \P_{\sqrt n\, x}(\tau>l)\Big)^{1/p},
\end{equation}
for some $C>0$ that does not depend on $n$, nor on $x$ or $R$. Now put $\lambda=1/p$.
\qed
\end{proofsect}

From now on we assume that a sufficiently high moment of the walk's steps is finite, in accordance with Proposition~\ref{Xtaukint}, see \eqref{muchoice}.

\begin{lemma} [Asymptotic relation between $V$ and $\Delta$]\label{lem-VasyDelta} Assume that the $\mu_k$-th moment of the steps is finite, with $\mu_k$ 
as in Proposition \ref{Xtaukint}. For any $M>0$, uniformly for $x\in W$ satisfying $|x|\leq M$,
\begin{equation}\label{VasyDelta}
\lim_{n\to\infty}n^{-\frac k4 (k-1)} V(\sqrt n\, x)=\Delta(x).
\end{equation}
\end{lemma}

\begin{proofsect}{Proof} For notational reasons, we write $m$ instead of $n$. As is seen from \eqref{Vdef}, it suffices to show that $\lim_{m\to\infty}m^{-\frac k4 (k-1)}\E_{\sqrt m\, x}[|\Delta(X(\tau))|]=0$, uniformly for $x\in W$ satisfying $|x|\leq M$. With some large $R>0$, we split this expectation into $\{\tau\leq Rm\}$ and $\{\tau> Rm\}$. For the first term, we are going to use Donsker's invariance principle. Indeed, under $\P_{\sqrt m\, x}$, the vector $B^{\ssup{m}}=(B_1^{\ssup{m}},\dots,B_k^{\ssup{m}})$ of the processes $B_i^{\ssup{m}}(t)=m^{-\frac 12} X_i(\lfloor tm\rfloor)$ converges towards a Brownian motion $B$ on $\R^k$ starting at $x$. Hence, for any $R>0$,  $m^{-\frac k4 (k-1)}\Delta(X(\tau))\1_{\{\tau\leq Rm\}}$ converges weakly and in $L^1$ towards $\Delta(B(T))\1_{\{T\leq R\}}$, which is equal to zero by continuity. It is not difficult to see that this convergence is uniform for $x\in W$ satisfying $|x|\leq M$.

The main difficulty is to show that $\limsup_{m\to\infty}m^{-\frac k4 (k-1)}\E_{\sqrt m\, x}[|\Delta(X(\tau))|\1_{\{\tau> Rm\}}]$ vanishes as $R\to\infty$, uniformly for $x\in W$ satisfying $|x|\leq M$. According to Lemma~\ref{integrneu}, we only have to show that, for $r\in(0,\frac k4(k-1)-1)$ and $\lambda\in(0,1)$ as in Lemma~\ref{integrneu},
\begin{equation}\label{integrneu1}
\limsup_{m\to\infty}m^{-\frac k4 (k-1)}\Big(\sum_{n\in\N}n^r \P_{\sqrt m\, x}(\tau>n)\Big)^\lambda<\infty,
\end{equation}
uniformly for $x\in W$ satisfying $|x|\leq M$. For this, is suffices to find, for any sufficiently small $\gamma>0$ (only depending on $r$ and $k$), some $C>0$ such that
\begin{equation}\label{tau>lesti}
\P_{\sqrt m\, x}(\tau>n)\leq C\Big(\frac {m^{1/\lambda}}n\Big)^{\frac k4(k-1)},\qquad \mbox{if }m=o(n^{1-\gamma}),
\end{equation}
uniformly for $x\in W$ satisfying $|x|\leq M$. 

The proof of \eqref{tau>lesti} uses parts of the proof of Proposition~\ref{Xtaukint}. We again restrict ourselves to the non-lattice case. We use $C$ as a generic positive constant that is uniform in $m$, $n$ and $x$ in the ranges considered. As in the proof of Proposition~\ref{Xtaukint}, we pick small parameters 
$\eta,\eps,\xi_1,\xi_2>0$ satisfying $8\eta<\xi_1$ and $\xi_2>\eps+\eta+\xi_1\frac k4(k-1)$. We assume that the $\mu$-th steps of the walk are finite, 
where $\mu$ is so large that \eqref{muchoice} holds.  
Abbreviate $a_{m,n}=\P_{\sqrt m\, x}(\tau>n)(n/m^{1/\lambda})^{\frac k4(k-1)}$ for $n,m\in\N$. For any $n\in\N$, we pick $m_n\leq o(n^{1-\gamma})$ maximal for $m\mapsto a_{m,n}$, and we put $A_n=\max\{a_{m_1,1},\dots,a_{m_n,n}\}$. Then our goal is to prove that $(A_n)_{n\in\N}$ is bounded. The index $m_n$, and hence also $A_n$, depends on $x$, but our estimates will be uniform in $x\in W$ satisfying $|x|\leq M$.

We split into the event $\{|X(n)|\leq n^{\frac 12+\eta}\}$ and the remainder. For the first term, we also use Lemma~\ref{KMGrephrase} with $f=\1$. This gives
\begin{equation}\label{tau>lesti1}
\begin{aligned}
\P_{\sqrt {m_n}\, x}(\tau>n)&\leq \P_{\sqrt {m_n}\, x}\big(|X(n)|>n^{\frac 12+\eta}\big)\\
&\quad +\int_W \d y\, \1_{\{|y|\leq n^\eta\}}\Big[ D_n^{\ssup{n}}(\sqrt {m_n}\, x, y)-\E_{\sqrt {m_n}\, x}\big[\1_{\{\tau \leq n\}}D_{n-\tau}^{\ssup{n}}(X(\tau), y)\big]\Big].
\end{aligned}
\end{equation}
Since in particular $m_n=o(n)$, the first term, $\P_{\sqrt {m_n}\, x}(|X(n)|>n^{\frac 12+\eta})$, can be estimated via Markov inequality against $C n^{-\mu \eta}$
using \eqref{momentnew}. If we choose $\mu$ such that $-\mu \eta + \frac k4 (k-1) <0$ (which is fulfilled under the moment condition in Proposition~\ref{Xtaukint}, see
\eqref{muchoice}), we obtain the bound  $\P_{\sqrt {m_n}\, x}(|X(n)|>n^{\frac 12+\eta}) \leq C n^{-\frac k4(k-1)} \leq C (m_n^{1/\lambda}/n)^{\frac k4(k-1)}$.
Let us turn to the second line of \eqref{tau>lesti1}.

As in the proof of Proposition~\ref{Xtaukint}, we split the expectation in the second line of \eqref{tau>lesti1} into the parts where $\tau\leq t_n$, $t_n\leq \tau\leq n-s_n$, and $n-s_n\leq \tau\leq n$, where $t_n=n^{1-\xi_1}$ and $s_n=n^{\frac 12+\xi_2}$. We want to apply Lemma~\ref{CLTtool} to the first part (together with the term $D_n^{\ssup{n}}(\sqrt {m_n}\, x, y)$) and Lemma~\ref{middleterm} to the second, i.e., we replace $x$ by $\sqrt {m_n} x$ in that lemmas. 
Lemma~\ref{middleterm} immediately applies since $\sqrt{m_n} \, x = o(\sqrt{n})$, and Lemma~\ref{CLTtool} applies if we assume 
that $\gamma>\xi_1$ to ensure that $\sqrt{m_n} \,|x| = o(\sqrt{t_n})$, which we do henceforth. 
Furthermore, for the last term we use \eqref{bdtotalmass} and obtain, as in \eqref{I-II-IIIesti}:
\begin{equation}\label{tau>lesti2}
\begin{aligned}
\P_{\sqrt {m_n}\, x}(\tau>n)&\leq Cn^{-\frac k4(k-1)}\int_W\1_{\{|y|\leq n^\eta\}}\frac{{\rm e}^{-\frac 12 |y|^2}}{(2\pi)^{k/2}}\Delta(y) \,\d y\,V_{t_n}\big(\sqrt {m_n} x\big)\\
&\qquad+ \Ocal\big(n^{\eps+\eta-\xi_2}\big)\P_{\sqrt{m_n}\, x}(\tau\geq t_n)+C\P_{\sqrt{m_n}\, x}(n-s_n\leq \tau\leq n) \\
&\qquad +o\big(n^{\frac 12-\frac\mu 4+\eta k+\xi_2}\big)+\Ocal(n^{1-\mu(\eta+{\xi_1}/4)}) +\Ocal\big(n^{1-\eps\mu}\big)+C\Big(\frac {m_n^{1/\lambda}}n\Big)^{\frac k4(k-1)}.
\end{aligned}
\end{equation}
Since $\mu$ satisfies \eqref{muchoice}, the last line is not larger than $C (m_n^{1/\lambda}/n)^{\frac k4(k-1)}$. Hence, we obtain
\begin{equation}\label{tau>lesti3}
\begin{aligned}
a_{{m_n},n}&\leq C m_n^{-\frac 1\lambda\,\frac k4(k-1)}\,V_{t_n}\big(\sqrt {m_n} x\big)+C n^{\eps+\eta-\xi_2}\Big(\frac n{t_n}\Big)^{\frac k4(k-1)}a_{m_n,t_n}\\
&+C \Big(a_{m_n,n-s_n} \Big(\frac n{n-s_n}\Big)^{\frac k4(k-1)}-a_{m_n,n}\Big)+C.
\end{aligned}
\end{equation}
Note that the factor $n^{\eps+\eta-\xi_2}(n/t_n)^{\frac k4(k-1)}=n^{\eps+\eta-\xi_2+\xi_1\frac k4(k-1)}$ is $o(1)$ by our requirement that the exponent is negative.  
We use Lemma~\ref{integr} to estimate, for 
$r\in(0,\frac k4(k-1)-1)$, $\lambda\in(0,1)$ and $a\geq 0$ as in that lemma,
$$
\begin{aligned}
V_{t_n}\big(\sqrt {m_n} x\big)&\leq \Delta\big(\sqrt {m_n} x\big)+C m_n^{(1+a)\frac k4(k-1)}+C\Big(\sum_{l=\lceil m_n^{1+a} \rceil}^{t_n} 
l^r\P_{\sqrt {m_n}\, x}(\tau>l)\Big)^\lambda\\
&\leq C m_n^{(1+a)\frac k4(k-1)}+C m_n^{(1+a) \frac k4(k-1)}\Big(\sum_{l= \lceil m_n^{1+a} \rceil}^{t_n} l^{r-\frac k4(k-1)}a_{m_n,l}\Big)^\lambda\\
&\leq C m_n^{(1+a)\frac k4(k-1)}+C m_n^{(1+a) \frac k4(k-1)}A_{n}^\lambda.
\end{aligned}
$$
Substituting this in \eqref{tau>lesti3} and solving for $a_{m_n,n}$, we obtain
\begin{equation}\label{tau>lesti4}
\begin{aligned}
a_{m_n,n}&\leq C+Cm_n^{(1+a-\frac 1\lambda)\frac k4(k-1)}+ C m_n^{(1+a)\lambda-\frac 1\lambda) \frac k4(k-1)}A_{n}^\lambda+o(1) a_{m_n,t_n}+\frac {C+o(1)}{C+1}a_{m_n,n-s_n}.
\end{aligned}
\end{equation}
Now pick $a>0$ small enough such that the second term on the right hand side is $\leq C$ and such that the third is $\leq A_{n}^\lambda$.
Recall that $m_n$ satisfies $a_{m_n,n}=\max_{m\leq o(n^{1-\gamma})}a_{m,n}$.
Picking $\xi_1>0$ even smaller, we can assume that $m_n =o( t_n^{1-\gamma})$, hence $a_{m_n,t_n}\leq a_{m_{t_n},t_n}\leq A_{t_n}\leq A_n$. Similarly, $a_{m_n,n-s_n}\leq A_n$. Hence, we obtain
$$
a_{m_n,n}\leq C+A_{n}^\lambda+\frac {C+o(1)}{C+1}A_n.
$$
Since the right hand side is increasing in $n$, we also obtain this estimate for $A_n$ instead of $a_{m_n,n}$. From this, it follows that $(A_n)_{n\in\N}$ 
is bounded. This finishes the proof.
\qed
\end{proofsect}

In the following lemma, we see in particular that $V$ does not increase much faster than $\Delta$ on $W$ at infinity. In particular, we can prove some integrability property of $V$, its regularity and its positivity. Recall that the $\mu$-th moment of the steps is assumed finite for some sufficiently large $\mu$, properly chosen in accordance with Proposition~\ref{Xtaukint}.

\begin{lemma}[Bounds on $V$, integrability, regularity and positivity]\label{Vbound} 
Assume that the $\mu_k$-th moment of the steps is finite, with $\mu_k$ as in Proposition \ref{Xtaukint}.
\begin{enumerate}
\item[(i)] There is a constant $C>0$ such that $V(x)\leq \Delta(x)+|x|^{\frac k2(k-1)}+C$ for any $x\in W$.

\item[(ii)] Fix $0<\nu\leq \frac\mu{\frac k2(k-1)}$, then $\E_x[V(X(n))^{\nu}\1_{\{\tau>n\}}]$ is finite for any $n\in\N$ and $x\in W$.

\item[(iii)] $V$ is regular for the restriction of the transition kernel to $W$.

\item[(iv)] $V$ is positive on $W$. 
\end{enumerate}
\end{lemma}

\begin{proofsect}{Proof} (i) According to Lemma~\ref{lem-VasyDelta}, there is $N_0\in\N$ such that, for any $n\in\N$ satisfying $n\geq N_0$ and for any $x\in W$ satisfying $|x|\leq 1$,
$$
V\big(x\sqrt n\big)\leq n^{\frac k4(k-1)}\big[\Delta(x)+1\big].
$$
Now let $x\in W$ be arbitrary. If $|x|\geq N_0+1$, then the above implies that 
$$
V(x)=V\Big(\frac x{\sqrt{\lceil |x|^2\rceil}}\,\sqrt{\lceil |x|^2\rceil}\Big)\leq \lceil |x|^2\rceil^{\frac k4(k-1)}\Big[\Delta\Big(\frac x{\sqrt{\lceil |x|^2\rceil}}\Big)+1\Big]
\leq \Delta(x)+(|x|+1)^{\frac k2(k-1)}.
$$
It suffices to show that $V$ is bounded on bounded subsets of $W$. It is clear that the map $x\mapsto \E_x[|\Delta(X(\tau))|\1_{\{\tau\leq 2\}}]$ is bounded on bounded  subsets of $W$. Use Lemma~\ref{integrneu} with $R=2$ and $n=1$ to estimate, for $x$ in some bounded subset of $W$,
$$
\E_x\big[|\Delta(X(\tau))|\1_{\{\tau> 2\}}\big]\leq C \E_x[\tau^{r+1}]^\lambda,
$$
see \eqref{schoen}. It is clear that the map $t\mapsto \E_{tx}[\tau^{r+1}]$ is increasing, since, for $t_1<t_2$, the random variable $\tau$ is stochastically smaller under $\P_{t_1x}$ than under $\P_{t_2 x}$. In the proof of Lemma~\ref{lem-VasyDelta} (see \eqref{integrneu1}) it is in particular shown that $x\mapsto \E_{tx}[\tau^{r+1}]$ is bounded on bounded subsets of $W$ if $t$ is sufficiently large. This ends the proof of (i).

(ii) By (i), we have, for any $\nu>0$,
$$
\E_x[V(X(n))^{\nu}\1_{\{\tau>n\}}]\leq \E_x\big[|\Delta(X(n))|^\nu\big]+\E_x\big[|X(n)|^{\nu\frac k2 (k-1)}\big]+C.
$$
Since $\Delta$ is a polynomial of degree $k-1$ in any $x_i$ and by independence of the components, the right hand side is finite as soon as both $\nu(k-1)$ and $\nu\frac k2 (k-1)$ do not exceed $\mu$. Since we assumed that $k\geq 2$, this is true as soon as $\nu\leq \frac\mu{\frac k2(k-1)}$.

(iii) We recall from \cite[Th.~2.1]{KOR02} that the process $(\Delta(X(n)))_{n\in\N_0}$ is a martingale under $\P_x$ for any $x\in \R^k$. In particular, $\E_x[\Delta(X(n))]=\Delta(x)$ for any $n\in\N$ and any $x\in\R^k$. The regularity of $V$ is shown as follows. For any $x\in W$,
\begin{equation}\label{harmonic}
\begin{aligned}
\E_x\bigl[\1_{\{\tau>1\}} V(X(1))\bigr]&=\E_x\bigl[\1_{\{\tau>1\}} \Delta(X(1))\bigr]-\E_x\bigl[\1_{\{\tau>1\}} \E_{X(1)}[\Delta(X(\tau))]\bigr]\\
&=\E_x\bigl[\1_{\{\tau>1\}} \Delta(X(1))\bigr]-\E_x\bigl[\1_{\{\tau>1\}}\Delta(X(\tau))\bigr]\\
&=\E_x\bigl[\1_{\{\tau>1\}} \Delta(X(1))\bigr]-\E_x[\Delta(X(\tau))]+\E_x\bigl[\Delta(X(\tau))\1_{\{\tau\leq 1\}}\bigr]\\
&=\E_x\bigl[\Delta(X(1))\bigr] - \E_x\bigl[\Delta(X(\tau))\bigr] + \E_x \bigl[\Delta(X(\tau)) \1_{\{\tau \leq 1\}} - \Delta(X(1)) \1_{\{\tau \leq 1\}}\bigr]  \\
&=V(x),
\end{aligned}
\end{equation}
where the second equality follows from the strong Markov property at time $\tau$.

(iv) Recall that $Y(n)=X(n)-X(n-1)\in\R^k$ is the step vector of the random walk at time $n$. Certainly, $Y(1)$ lies in $\overline W$ with positive probability. From the Centering Assumption it follows that $Y_1(1)$ is not constant almost surely. Therefore, the vector $v=\E[Y(1)\mid Y(1)\in\overline W]$ lies in $W$, since, for any $i=1,\dots,k-1$, we have $v_i-v_{i-1}=\E[Y_i(1)-Y_{i-1}(1)\mid Y_i(1)-Y_{i-1}(1)\geq 0]$, and $Y_i(1)-Y_{i-1}(1)$ is positive with positive probability on $\{Y_i(1)-Y_{i-1}(1)\geq 0\}$. 

Let $A$ be a closed neighborhood of $v$ that is contained in $W$. Hence, for any sufficiently large $m \in \N$, we have that 
$$
\begin{aligned}
\P_x(\tau > m,X(m)\in m A)&\geq \P_x\big(Y(1),\dots,Y(m)\in \overline W\big)\P_x\Big(\frac 1m X(m)\in A\,\Big|\,Y(1),\dots,Y(m)\in \overline W\Big)>0,
\end{aligned}
$$
since the first term is positive for any $m$, and the last one converges to one, according to the weak law of large numbers. According to Lemma~\ref{lem-VasyDelta}, for any sufficiently large $n\in\N$ and for any $y\in A$, $V(\sqrt n y)\geq \frac 12 n^{-\frac k4(k-1)}\inf_A\Delta$. In particular, $\inf_{y\in m A}V(y)>0$ for any sufficiently large $m$.

Now recall from Corollary~\ref{Vpos} that $V\geq 0$ and iterate the regularity equation for $V$ to the effect that 
$$
\begin{aligned}
V(x)&=\E_x[V(X(m))\1_{\{\tau>m\}}]\geq  \E_x[V(X(m))\1_{\{X(m)\in m A\}}\1_{\{\tau>m\}}]\\
&\geq \inf_{m A}V\,\,\P_x(X(m)\in m A, \tau>m)\\
&>0.
\end{aligned}
$$
Hence, $V(x)$ is positive.
\qed
\end{proofsect}

\begin{bem}[$(V_n)_{n\in\N_0}$ as an iterating sequence]
A modification of the calculation in \eqref{harmonic} shows that $\E_x\bigl[\1_{\{\tau>1\}} V_n(X(1))\bigr]=V_{n+1}(x)$ for any $x\in W$ and $n\in \N$. Furthermore, it is clear that $V_0=\Delta$. In other words, we can see the sequence $(V_n)_{n\in\N_0}$ as the iterating sequence for the iterated application of the expectation before the first violation of the ordering, starting with initial function $\Delta$.
\hfill$\Diamond$
\end{bem}

Now that we know that $V\colon W\to(0,\infty)$ is a positive regular function for the restriction of the transition kernel to $W$, we can finally define the Doob $h$-transform of $X$ on $W\cap S^k$ with $h=V$. Recalling \eqref{Dndef}, its transition probabilities are given by
\begin{equation}\label{transfProc}
\begin{aligned}
\widehat\P^{\ssup V}_x(X(n)\in\d y)&=\P_x(\tau>n;X(n)\in\d y)\frac{V(y)}{V(x)}\\
&=\Bigl[\Dcal_n(x,\d y)-\E_x\bigl[\1_{\{\tau\leq n\}}\Dcal_{n-\tau}(X(\tau),\d y)\bigr]\Bigr]\frac{V(y)}{V(x)},\qquad n\in\N, x,y\in W.
\end{aligned}
\end{equation}

The measure on the right hand side is indeed a probability measure in $\d y$ on $W$, since, by Lemma~\ref{Vbound}(ii) it has finite mass on $W$, and by 
Lemma~\ref{Vbound}(iii) its mass is even equal to one.

Now we can show that the transformed process deserves the name \lq $k$ random walks conditioned on being strictly ordered for ever\rq.

\begin{lemma}[Conditional interpretation]\label{lem-CondInt} Assume that the $\mu_k$-th moment of the steps is finite, with $\mu_k$ as in Proposition \ref{Xtaukint}. 
The conditional distribution of the process $(X(n))_{n\in\N_0}$ given $\{\tau>m\}$ converges, as $m\to\infty$, to the Doob $h$-transform of $(X(n))_{n\in\N_0}$ with $h=V$, i.e., for any $x\in W$ and $n\in\N$,
\begin{equation}
\lim_{m\to\infty}\P_x(X(n)\in\d y\mid\tau>m)=\widehat\P^{\ssup V}_x(X(n)\in\d y),\qquad \mbox{weakly.}
\end{equation}
\end{lemma}

\begin{proofsect}{Proof} Using the definition of the conditional probability and the Markov property at time $n$, we see that, for any $n,m\in\N$ satisfying $n<m$,
$$
\P_x(X(n)\in\d y\mid\tau>m)=\frac{\P_x(\tau>n;X(n)\in\d y)m^{\frac k4(k-1)}\P_y(\tau>m-n)}{m^{\frac k4(k-1)}\P_x(\tau>m)}.
$$
According to \eqref{tauasy}, the last term in the numerator converges towards $KV(y)$, and the denominator converges towards $KV(x)$ as $m \to \infty$. 
Compare to the first line of \eqref{transfProc} to see that this finishes the proof.
\end{proofsect}
\qed

Recall that Dyson's Brownian motions is the Doob $h$-transform of a standard Brownian motion on $W$ with $h$ equal to the restriction of the 
Vandermonde determinant $\Delta$ to $W$. Recall that $\widehat{\P}_x^{\ssup V}$ 
is the Doob $h$-transform with $h=V$ of the random walk $X$ on $W\cap S^k$. 
The Brownian motion $B=(B_1,\dots,B_k)$ on $\R^k$ starts from $x\in W$ under ${\tt P}_x$, and ${\tt E}_x$ denotes the corresponding expectation.
Denote by ${\tt P}_x^{\ssup \Delta}(B(t)\in\d y)$ the Doob $h$-transform with $h=\Delta$ of the Brownian motion $B$ on $\R^k$, e.g.
$$
{\tt P}_x^{\ssup \Delta}(B(t)\in\d y) = \P_x(T>t;B(t)\in\d y)\frac{\Delta(y)}{\Delta(x)}, \quad x \in W, t>0. 
$$

\begin{lemma}[Convergence towards Dyson's Brownian motions]\label{lem-Dysonconv} Assume that the $\mu_k$-th moment of the steps is finite, for $\mu_k$ as in 
Proposition~\ref{Xtaukint}. Then, under $\widehat{\P}_{x\sqrt n}^{\ssup V}$, the process $B^{\ssup n}=(n^{-\frac 12}X(\lfloor tn\rfloor))_{t\in[0,\infty)}$ weakly converges, as $n\to\infty$, towards Dyson's Brownian motions started at $x$. More precisely, the sequence $(B^{\ssup n})_{n\in\N}$ is tight, and, for any $x\in W$, and any $t>0$,
\begin{equation}\label{Dysonconv}
\lim_{n\to\infty} \widehat{\P}_{x\sqrt n}^{\ssup V}\Big(\frac 1{\sqrt n}X(\lfloor tn\rfloor)\in \d y\Big)
={\tt P}_x^{\ssup \Delta}(B(t)\in\d y),\qquad \mbox{weakly.}
\end{equation}
\end{lemma}

\begin{proofsect}{Proof}
Using \eqref{transfProc} and Lemmas~\ref{lem-Invariance} and \ref{lem-VasyDelta}, we see that, for any $t>0$, as $n\to\infty$,
\begin{equation}
\begin{aligned}
\widehat{\P}_{x\sqrt n}^{\ssup V}\Big(\frac 1{\sqrt n}X(\lfloor tn\rfloor)\in \d y\Big)
&=\P_{x\sqrt n}\Big(\tau>\lfloor tn\rfloor;\frac 1{\sqrt n}X(\lfloor tn\rfloor)\in \d y\Big)\frac{V(\sqrt n \,y)}{V(\sqrt n \,x)}\\
&\to \P_x(T>t;B(t)\in\d y)\frac{\Delta(y)}{\Delta(x)}={\tt P}_x^{\ssup \Delta}(B(t)\in\d y).
\end{aligned}
\end{equation}
This shows that \eqref{Dysonconv} holds.

Now we show the tightness. According to the Kolmogorov-Chentsov criterion, it suffices to find, for any $S>0$, constants $\alpha,\beta, C>0$ such that 
\begin{equation}\label{tight}
{\E}_{x\sqrt n}^{\ssup V}\big[\big|B^{\ssup n}(t)-B^{\ssup n}(s)\big|^\alpha\big]\leq C |t-s|^{1+\beta},\qquad s,t\in[0,S], n\in\N.
\end{equation}
This is done as follows. We pick some $\alpha\in(2,4)$. Fix $0\leq s<t\leq S$. First note that using the Markov property we obtain
\begin{equation}\label{tight1}
\begin{aligned}
{\E}_{x\sqrt n}^{\ssup V}&\big[\big|B^{\ssup n}(t)-B^{\ssup n}(s)\big|^\alpha\big]\\
&=\int_W\int_W |z_1-z_2|^\alpha\,\P_{x\sqrt n}\Big(\tau>\lfloor sn\rfloor ,\frac{X(\lfloor sn\rfloor)}{\sqrt n}\in\d z_1\Big)\\
&\qquad\qquad\times\P_{z_1\sqrt n}\Big(\tau>\lfloor tn\rfloor-\lfloor sn\rfloor ,\frac{X(\lfloor tn\rfloor-\lfloor sn\rfloor)}{\sqrt n}\in\d z_2\Big)\frac{V(z_2\sqrt n)}{V(x\sqrt n)}\\
&\leq \int_W\int_W |z_1-z_2|^\alpha\,\P_{x\sqrt n}\Big(\frac{X(\lfloor sn\rfloor)}{\sqrt n}\in\d z_1\Big)\P_{z_1\sqrt n}\Big(\frac{X(\lfloor tn\rfloor-\lfloor sn\rfloor)}{\sqrt n}\in\d z_2\Big)\frac{V(z_2\sqrt n)}{V(x\sqrt n)}.
\end{aligned}
\end{equation}
We use $C$ as a generic positive constant, not depending on $s,t$ (as long as $0\leq s<t\leq S$) nor on $n$, nor on $z_1$ or $z_2$.
By Lemma~\ref{lem-VasyDelta}, $1/V(x\sqrt n)\leq C n^{-\frac k4(k-1)}$, uniformly in $x$ on compact subsets of $W$. From Lemma~\ref{Vbound}(i), we know that there is a polynomial $P\colon\R^k\to\R$ of degree $\leq \frac k2(k-1)$ such that $V(z_2\sqrt n)\leq |P(z_2)| n^{\frac k4(k-1)}$ for any $n\in\N$ and $z_2\in W$. Hence,
$$
{\E}_{x\sqrt n}^{\ssup V}\big[\big|B^{\ssup n}(t)-B^{\ssup n}(s)\big|^\alpha\big]
\leq C{\E}_{x\sqrt n}\big[\big|B^{\ssup n}(t)-B^{\ssup n}(s)\big|^\alpha |P(B^{\ssup n}(t))|\big].
$$
Now we use H\"older's inequality with $p=4/\alpha$ and $q=4/(4-\alpha)$, to obtain
$$
{\E}_{x\sqrt n}^{\ssup V}\big[\big|B^{\ssup n}(t)-B^{\ssup n}(s)\big|^\alpha\big]
\leq C {\E}_{x\sqrt n}\big[\big|B^{\ssup n}(t)-B^{\ssup n}(s)\big|^4\big]^{\alpha/4} {\E}_{x\sqrt n}\big[|P(B^{\ssup n}(t))|^{4/(4-\alpha)}\big]^{1-\alpha/4}.
$$
It is known that the first expectation ${\E}_{x\sqrt n}\big[\big|B^{\ssup n}(t)-B^{\ssup n}(s)\big|^4\big]$ on the right hand side 
can be estimated against $C|t-s|^2$. Furthermore, the second expectation is bounded in $n\in\N$ and $t\in[0,1]$ as soon as the $(\frac k2(k-1)\frac 4{4-\alpha})$-th moment of the steps is finite, i.e., as soon as $\frac k2(k-1)\frac 4{4-\alpha}\leq \mu_k$. Choosing $\alpha$ sufficiently close to $2$, this is satisfied, by our assumption that $\mu_k>k(k-1)$. For this choice of $\alpha$ we obtain
$$
{\E}_{x\sqrt n}^{\ssup V}\big[\big|B^{\ssup n}(t)-B^{\ssup n}(s)\big|^\alpha\big]
\leq C  |t-s|^{\alpha/2},
$$
which shows that \eqref{tight} holds with $\beta=\alpha/2-1>0$.
\qed
\end{proofsect}

\end{document}